\documentclass[11pt,reqno]{amsart}
\usepackage{amssymb}
\usepackage{amsfonts}
\usepackage{amsmath}
\usepackage{graphicx}
\usepackage{hyperref}
\usepackage{float}
\usepackage{epstopdf}
\usepackage{color}
\usepackage{comment}
\usepackage{bm}
\usepackage{soul}

\allowdisplaybreaks
\newtheorem{theorem}{Theorem}[section]
\newtheorem{proposition}[theorem]{Proposition}
\newtheorem{lemma}[theorem]{Lemma}

\newtheorem{definition}[theorem]{Definition}

\newtheorem{thm}{Theorem}

\def\SG{\mathcal{SG}}

\numberwithin{equation}{section}

\begin{document}
\title[Sobolev spaces on p.c.f. self-similar sets: Boundary behavior and interpolation theorems]{Sobolev spaces on p.c.f. self-similar sets: Boundary behavior and interpolation theorems}

\author{Shiping Cao}
\address{Department of Mathematics, Cornell University, Ithaca 14853, USA}
\email{sc2873@cornell.edu}
\thanks{}

\author{Hua Qiu}
\address{Department of Mathematics, Nanjing University, Nanjing 210093, China}
\email{huaqiu@nju.edu.cn}
\thanks{The research of Qiu was supported by the NSFC grant 11471157}

%    General info

\date{}


\begin{abstract}
We study the Sobolev spaces $H^\sigma(K)$ and $H^\sigma_0(K)$ on p.c.f. self-similar sets in terms of the boundary behavior of functions. First, for $\sigma\in \mathbb{R}^+$, we make an exact description of the tangents of functions in $H^\sigma(K)$ at the boundary.   Second, we characterize $H_0^\sigma(K)$ as the space of functions in $H^\sigma(K)$ with zero tangent of an appropriate order depending on $\sigma$. Last, we extend  $H^\sigma(K)$ to $\sigma\in\mathbb{R}$, and obtain various interpolation theorems with $\sigma\in\mathbb{R}^+$ or $\sigma\in\mathbb{R}$. We illustrate that there is a countable set of critical orders, that arises naturally in the boundary behavior of functions, such that $H^\sigma_0(K)$ presents a critical phenomenon if $\sigma$ is critical. These orders will play a crucial role in our study. They are just the values in $\frac 12+\mathbb{Z}_+$ in the classical case, but are much more complicated in the fractal case.
\end{abstract}
\maketitle

\tableofcontents

\section{Introduction}\label{intro}

The boundary behavior of functions, as an important topic in analysis on fractals, has been studied for years since the construction of the Laplacians on fractals.  See \cite{k2,k3} for Kigami's construction of the Laplacians on p.c.f. self-similar sets, see \cite{BB1,BB2,BBKT,G,KZ,Lindstrom} for the probabilistic approach, and see books \cite{kigami,s2} for further developments. Many important results are obtained, including tangents and gradients \cite{k1,K1,calculus,distribution,s3,T1,cq1,cq4,DRS}. See also  \cite{smoothbump} and \cite{distribution} for related topics on the smooth bump functions and distributions on fractals.

In this paper, we will take a further step to study the boundary behavior of functions in Sobolev spaces on p.c.f. self-similar sets, which are analogs of $H^\sigma(\Omega), \sigma\in\mathbb{R},$ in the $\mathbb{R}^n$ case.  This is a continuation of our previous work \cite{cq3}.

 Recall that for a domain $\Omega$ in $\mathbb{R}^n$ with smooth boundary, for $\sigma\in\mathbb{Z}_+$, $H^\sigma(\Omega)$ is the space of function $f$'s on $\Omega$ such that $f$ and its derivatives (in the sense of distributions) up to order $\sigma$ are in $L^2(\Omega)$, and the definition of $H^\sigma(\Omega)$ can be generalized to all real $\sigma$ {via complex interpolation or other numerous equivalent methods.} There are rich fundamental results concerning the boundary behavior of functions in $H^\sigma(\Omega)$, including the trace theorems, interpolation theory, which provide powerful tools in the study of non-homogeneous boundary value problems and further topics on $\Omega$. Among them, the characterization of $H_0^\sigma(\Omega)$,  $\sigma\geq 0$,
 \begin{equation}\label{h0char}H^\sigma_0(\Omega)=\big\{f\in H^\sigma(\Omega): \frac{\partial ^j f}{\partial {{\nu}}^j}=0, \quad \forall 0\leq j<\sigma-\frac 12\big\},\end{equation}
 which was first discovered in \cite{LM2,LM3,LM4} by J.L. Lions and E. Magenes,  plays a central and delicate role.  See monograph \cite{Lion} for a systematic development and various applications. 

We will reproduce the characterization (\ref{h0char}) of $H_0^\sigma(\Omega)$ in the fractal setting, and stem from which,  we aim to provide a throughout study on the interpolation of $H^\sigma(\Omega)$ on fractals. Due to the complicity of the fractal feature, we need to take a quite different approach. 

In the classical setting,  for $\sigma\notin \frac{1}{2}+\mathbb{Z}_+$, $H_0^\sigma(\mathbb{R}_+^n)$ can be embedded into $H^\sigma(\mathbb{R}^n)$ as a subspace, consisting of functions with support in $\mathbb{R}_+^n$, by extending functions by $0$ outside $\mathbb{R}_+^n$. What's more, there exists a retraction mapping $T:H^\sigma(\mathbb{R}^n)\to H_0^\sigma(\mathbb{R}_+^n)$. The proof of (\ref{h0char}) and the interpolation result of $H^\sigma(\Omega)$ essentially rely on this extension, and the local coordinate representation of $H^\sigma(\Omega)$ along the boundary. The values in $\frac 12+\mathbb{Z}_+$ are called critical orders, since $H_0^\sigma(\Omega)$ will present some critical phenomena when $\sigma$ is such a value.

However, on the p.c.f. self-similar sets, there are `derivatives' other than Laplacians and normal derivatives at the boundary, Although these new derivatives do not matter in the matching conditions when {extending} a function to a larger fractal domain, they still reflect the boundary behavior of functions \cite{s3}. As a consequence, for a p.c.f. fractal domain $\Omega$ with boundary, $H^\sigma_0(\Omega)$ no longer has the nice characterization as the space of functions on a larger fractal domain with support in $\Omega$, and the retraction mapping does not exist. In addition, the occurrence of the new derivatives will create more `critical orders' in the fractal setting.

Instead of extending functions by $0$ outside the fractal, we will {extract the information of their boundary behavior with a more straightforward method, by splitting the functions.} See Theorem \ref{thm47}, Theorem \ref{thm410} and the remark after Definition \ref{def55}.  This method features our work, and shows natural insight into the Sobolev spaces. Here we summary the decomposition of spaces by splitting, but leave the explanation of notations in later context.  

\begin{thm}\label{thm10}
	Let $K$ be a p.c.f. self-similar set with boundary $V_0$. Let $k\geq 1$ be an integer, and $0\leq\sigma\leq 2k$. We have 
	 \[H^\sigma(K\setminus V_0)=\ker_\sigma\mathcal{T}_{\mathcal{H}_{k-1}}\oplus\big(\oplus_{\omega\in \mathcal{P}}\mathcal{R}_{\mathcal{H}_{k-1},\omega}\big(l^2(\mathcal{H}_{k-1},A_w;r_w^{\sigma/2}\mu_w^{(\sigma-1)/2})\big)\big),\]
\[H^\sigma_0(K\setminus V_0)=\ker_\sigma\mathcal{T}_{\mathcal{H}_{k-1}}\oplus\big(\oplus_{\omega\in \mathcal{P}}\mathcal{R}_{\mathcal{H}_{k-1},\omega}\big(\overline{l^2(\mathcal{H}_{k-1};r_w^{\sigma/2}\mu_w^{(\sigma-1)/2})}\big)\big),\]
	\[H^\sigma_{00}(K\setminus V_0)=\ker_\sigma\mathcal{T}_{\mathcal{H}_{k-1}}\oplus\big(\oplus_{\omega\in \mathcal{P}}\mathcal{R}_{\mathcal{H}_{k-1},\omega}\big(l^2(\mathcal{H}_{k-1};r_w^{\sigma/2}\mu_w^{(\sigma-1)/2})\big)\big).\]
\end{thm}

A surprising consequence of the above decomposition result is that it  provides privilege when {considering} the interpolation couple $(H^{\sigma_1}_0(K\setminus V_0),H^{\sigma_2}_0(K\setminus V_0))$ for at least one of $\sigma_1,\sigma_2$ being a critical order, while in the $\mathbb{R}^n$ case, Lions and Magenes's method will meet teratological difficulty, see \cite{Lion} (Chapter 1, Section 18). For example, when $K$ is chosen to be the unit interval $I=[0,1]$, we will have no difficulty to generalize the interpolation result for $[H^{\sigma_1}(0,1),H^{\sigma_2}(0,1)]_\theta$ when $\sigma_1$ or $\sigma_2$ is in $\frac 12+\mathbb{Z}_+$. 

Now we briefly introduce our main results. Let $K$ be a p.c.f. self-similar set which possesses a local regular Dirichlet form in the sense of Kigami. Let $V_0$ be its boundary consisting of finite points. For $\sigma\in\mathbb{R}$,  by a slight abuse of notation, we write $H^\sigma(K)$ for the Sobolev space on the domain $K\setminus V_0$. A systematical introduction of Sobolev spaces can be found in \cite{s1} on fractals  by R.S. Strichartz and in \cite{Gr} on more general metric measure spaces by A. Grigor'yan. See also \cite{cq3, GL, HZ} for some equivalent Besov type characterizations of $H^\sigma(K)$,  and \cite{CG1,CG, cq3, GKS} for related interpolation results. In this paper, we focus on the following three aspects.

Firstly, we study tangents at boundary points for functions in $H^\sigma(K)$ with $\sigma\in \mathbb{R}^+$. In history, various different approaches are developed towards gradients and tangents for functions on $K$. Typical ideas include defining the gradients by the energy measures \cite{k1,K1} and defining the tangents as the multiharmonic functions that match the local behavior of functions $f$ at a generic point \cite{T1} or a vertex \cite{s3,distribution} in $K$. We will introduce a simpler but more efficient definition based on {the latter} idea, and give a thorough study of tangents at points in $V_0$ for functions in $H^\sigma(K)$. See Definition \ref{def33} and Theorem \ref{thm314}.  

Secondly, we study the Sobolev space $H^\sigma_0(K)$ with $\sigma\in\mathbb{R}^+$, which is defined as the closure of all compactly supported smooth functions with respect to the norm of $H^\sigma(K)$. In particular, we will show that(Theorem \ref{thm42}), analogously to (\ref{h0char}),
\begin{thm}\label{thm11}
For $\sigma\geq 0$, $H^\sigma_0(K)=\big\{f\in H^\sigma(K):T^{(\sigma)}_\omega(f)=0,\forall \omega\in \pi^{-1}(V_0)\big\}.$ In particular, $H^\sigma_0(K)=H^\sigma(K)$ if $\sigma\leq \frac{d_S}{2}$.
\end{thm}
Here $\pi$ is the canonical coding map associated with $K$, $T^{(\sigma)}_\omega(f)$(Definition \ref{def33} and \ref{def313}) stands for the tangent of $f$ at $\pi(\omega)$ with order that, roughly speaking, works best for $H^\sigma(K)$, and $d_S$ is the spectral dimension of $K$. Readers are suggested to  compare this result with the authors' previous work on the characterization of $H^\sigma_D(K)$ and $H^\sigma_N(K)$ in \cite{cq3}. As we have mentioned, the proof of Theorem \ref{thm11} essentially relies on the splitting method in Theorem \ref{thm10}, and is very different from Lions and Magenes's method for classical domains in \cite{Lion}. We will use the smooth bump functions developed by L.G. Rogers, R.S. Strichartz and A. Teplyeav in \cite{smoothbump} as an important tool. 

Lastly, we study the interpolation theorems concerning $H^\sigma(K)$ and $H^\sigma_0(K)$. First, using the results obtained in the previous parts, we are ready to deal with the case $\sigma\in\mathbb{R}^+$(Theorem \ref{thm54},\ref{thm56} and \ref{thm57}),
\begin{thm}\label{thm12} 
	\[[H^\sigma(K),H^{\sigma'}(K)]_\theta=H^{(1-\theta)\sigma+\theta\sigma'}(K), \quad \forall \sigma>\sigma'\geq 0, \theta\in[0,1],\]
	\[[H^\sigma_0(K),H^{\sigma'}_0(K)]_\theta=H^{(1-\theta)\sigma+\theta\sigma'}_{00}(K), \quad \forall \sigma>\sigma'\geq 0, \theta\in(0,1),\] 
	where $H^\sigma_{00}(K)\subset H^\sigma_0(K)$ are analogs of the Lions-Magenes spaces. 
	In particular, 
	\[[H_{00}^\sigma(K),H_{00}^{\sigma'}(K)]_\theta=H_{00}^{(1-\theta)\sigma+\theta\sigma'}(K), \quad \forall \sigma>\sigma'\geq 0, \theta\in[0,1],\]
and $H^\sigma_{00}(K)=H^\sigma_0(K)$ except a countable set of critical orders of $\sigma$ that arises naturally in Theorem \ref{thm314} dealing with tangents of functions in $H^\sigma(K)$.
\end{thm}
 Moreover, we then introduce the space $H^\sigma(K)$ with $\sigma<0$ as the dual of $H^{-\sigma}_0(K)$, and extend the story of interpolation theorem to $\sigma\in\mathbb{R}$. The difficulty in this part lies in the fact that the domain of Laplacians is not closed under multiplication \cite{BST}, and we develop a projection technique that preserves regularity instead. It holds that(Theorem \ref{thm62}), 
\begin{thm}\label{thm13}
 For $-\infty<\sigma'<\sigma<\infty$ and $0<\theta<1$, 
	\[[H^{\sigma}(K),H^{\sigma'}(K)]_{\theta}=\begin{cases}
H^{\sigma_\theta}(K),&\text{ if }\sigma_\theta=(1-\theta)\sigma+\theta\sigma'\geq 0,\\
\big(H_{00}^{-\sigma_\theta}(K)\big)',&\text{ if }\sigma_\theta=(1-\theta)\sigma+\theta\sigma'<0.
\end{cases}\] 
In particular, $[H^{\sigma}(K),H^{\sigma'}(K)]_{\theta}=H^{\sigma_\theta}(K)$
except $\sigma_\theta$ is in the countable set of critical orders. 
\end{thm}

We briefly introduce the structure of our writing. In Section 2, we provide backgrounds and definitions that will be used later. {The} assumption (\textbf{A1}) will be introduced for convenience. In Section 3, we study the tangents at the boundary points for functions in  $H^\sigma(K)$ with $\sigma\in\mathbb{R}^+$. In Section 4, we characterize $H^\sigma_0(K)$  in terms of the boundary behavior of functions in $H^\sigma(K)$. In Section 5, we develop interpolation theorems for Sobolev spaces $H^\sigma(K),H^\sigma_0(K)$ and $H^\sigma_{00}(K)$ with $\sigma\in\mathbb{R}^+$. In Section 6, we extend the interpolation theorem of $H^\sigma(K)$ to $\sigma\in\mathbb{R}$. In Section 7, we present some examples, along with some equivalent narrations of our results. In the last section, we give a counter example of  (\textbf{A1}), and a brief discussion on how to prove the previous main theorems without assuming (\textbf{A1}).

\section{Preliminaries}
We introduce some backgrounds in this section, including the Dirichlet forms and Sobolev spaces on p.c.f. self-similar sets.

Let $\{F_i\}_{i=1}^N$ be a finite collection of contractions on a complete metric space $(X,d)$. The self-similar set associated with the \textit{iterated function system (i.f.s.)} $\{F_i\}_{i=1}^N$ is the unique compact set $K\subset X$ satisfying
$K=\bigcup_{i=1}^N F_iK.$
For $m\geq 1$, we define $W_m=\{1,\cdots,N\}^m$ the collection of \textit{words} of length $m$, and for each $w=w_1w_2\cdots w_m\in W_m$, denote
\[F_w=F_{w_1}\circ F_{w_2}\circ\cdots \circ F_{w_m}.\]
For uniformity, we set $W_0=\{\emptyset\}$, with $F_\emptyset$ being the identity map. For convenience, let $W_*=\bigcup_{m=0}^\infty W_m$ be the collection of all finite words.

Let  $\Sigma=\{1,2,\cdots,N\}^{\mathbb{N}}$ be the shift space endowed with the natural product topology. There is a continuous surjection $\pi: \Sigma\rightarrow K$ defined by
\[\pi(\omega)=\bigcap_{m\geq 1}F_{[\omega]_m}K,\]
where for $\omega=\omega_1\omega_2\cdots$ in $\Sigma$ we write $[\omega]_m=\omega_1\omega_2\cdots \omega_m\in W_m$ for each $m\geq 1$. Let 
\[C_K=\bigcup_{i\neq j}F_iK\cap F_jK,\quad \mathcal{C}=\pi^{-1}(C_K),\quad \mathcal{P}=\bigcup_{m\geq 1}\sigma^m \mathcal{C},\]
where $\sigma$ is the shift map define as $\sigma(\omega_1\omega_2\cdots)=\omega_2\omega_3\cdots$. $\mathcal{P}$ is called the \textit{post critical set}. Call $K$  a \textit{p.c.f. self-similar set} if $\#\mathcal{P}<\infty$. In what follows, we always assume that $K$ is a connected p.c.f.  self-similar set. 

Let $V_0=\pi(\mathcal{P})$ and call it the \textit{boundary} of $K$. For $m\geq 1$, we always have $F_w K\cap F_{w'}K\subset F_w V_0\cap F_{w'}V_0$ for any $w\neq w'\in W_m$. For simplicity, we assume (\textbf{A1}) throughout Section 4 to  6. \vspace{0.1cm}

\noindent(\textbf{A1}): \textit{For any $p\in V_0$, we assume $\#\pi^{-1}(p)=1$.}\vspace{0.1cm}

The assumption (\textbf{A1}) is a geometric assumption that provides some convenience in the following context, but it is not necessary. It is naturally satisfied for nested fractals, see \cite{kigami, Lindstrom}. In the next section, we will introduce another assumption (\textbf{A2}) on the measure $\mu$ defined on $K$. We will see in Section 8 that all our results are valid even if we do not assume (\textbf{A1}), as long as (\textbf{A2}) is assumed. See an example that (\textbf{A1}) fails in Section 8.

\subsection{Dirichlet forms on p.c.f. self-similar sets}
Let's briefly recall the construction of Dirichlet forms on p.c.f. self-similar sets. Readers are suggested to refer to  books \cite{kigami} and \cite{s2} for any unexplained details and notations. 

For $m\geq 1$, denote $V_m=\bigcup_{w\in W_m}F_wV_0$ and let $l(V_m)=\{f: f \text{ maps } V_m \text{ into } \mathbb{R}\}$. Write $V_*=\bigcup_{m\geq 0}V_m$.

Let $H=(H_{pq})_{p,q\in V_0}$ be a symmetric linear operator (matrix) on $l(V_0)$. $H$ is called a \textit{(discrete) Laplacian} on $V_0$ if  $H$ is non-positive definite; $Hu=0$ if and only if $u$ is constant on $V_0$; and $H_{pq}\geq 0$ for any $p\neq q\in V_0$.
Given a Laplacian $H$ on $V_0$ and a vector $\bm{r}=\{r_i\}_{i=1}^N$ with $r_i>0$, $1\leq i\leq N$, define the \textit{(discrete) Dirichlet form} on $V_0$ by
$$\mathcal{E}_0(f,g)=-(f,Hg),$$
for $f,g\in l(V_0)$, and inductively {{on $V_m$ by}}
$$\mathcal{E}_m(f,g)=\sum_{i=1}^Nr^{-1}_i\mathcal{E}_{m-1}(f\circ F_i, g\circ F_i), m\geq 1,$$
for $f, g\in l(V_m)$. Write {$\mathcal{E}_m(f):=\mathcal{E}_m(f,f)$} for short.

Say $(H,\bm{r})$ is a \textit{harmonic structure} if for any $f\in l(V_0)$,
\[\mathcal{E}_0(f)=\min\{\mathcal{E}_1(g):g\in l(V_1),g|_{V_{0}}=f\}.\]
{In addition, call $(H,\bm{r})$ a \textit{regular harmonic structure}, if $0<r_i<1, \forall 1\leq i\leq N$. In this paper, we will always assume that there exists a regular harmonic structure associated with $K$. }

{Now for each $f\in C(K)$, the sequence $\{\mathcal{E}_m(f)\}_{m\geq 0}$ is nondecreasing, so the following definitions make sense. Let $dom\mathcal{E}=\{f\in C(K):\lim\limits_{m\to\infty}\mathcal{E}_m(f)<\infty\},$ and 
\[
\mathcal{E}(f,g)=\lim_{m\to\infty} \mathcal{E}_m(f,g) \text{ for }f,g\in dom\mathcal{E}.
\]
We write $\mathcal{E}(f):=\mathcal{E}(f,f)$ for short, and call $\mathcal{E}(f)$ the \textit{energy} of $f$.
It is known that $(\mathcal{E},dom\mathcal{E})$ turns out to be a \textit{local regular Dirichlet form} on $L^2(K,\mu)$ for any  Radon measure $\mu$ on $K$.} 

An important feature of the form $(\mathcal{E},dom\mathcal{E})$ is the \textit{self-similar identity},
\begin{equation}\label{eq21}
\mathcal{E}(f,g)=\sum_{i=1}^Nr_i^{-1}\mathcal{E}(f\circ F_i, g\circ F_i), \quad \forall f,g\in dom\mathcal{E}.
\end{equation}
Denote $r_w=r_{w_1}r_{w_2}\cdots r_{w_m}$ for each $w\in W_m, m\geq 0$. Then for  $m\geq 1$, we have 
\begin{equation*}
\mathcal{E}_m(f,g)=\sum_{w\in W_m} r_w^{-1}\mathcal{E}_0(f\circ F_w, g\circ F_w),\quad \mathcal{E}(f,g)=\sum_{w\in W_m} r_w^{-1}\mathcal{E}(f\circ F_w, g\circ F_w).
\end{equation*}

Lastly, we need to mention that there is a natural metric on $K$ related with the energy form $(\mathcal{E},dom\mathcal{E})$, called the \textit{effective resistance metric}, which is defined as $$R(x,y)=\big(\min\{\mathcal{E}(f):f\in dom\mathcal{E} \text{ and }f(x)=1,f(y=0)\}\big)^{-1}, \forall x,y\in K.$$

\subsection{The Laplacians and Sobolev spaces}
Now we come to the basic concepts of the Laplacians and Sobolev spaces on $K$. Readers may find detailed backgrounds and further discussions in various contexts, for example \cite{cq3,kigami,s1,s2}. 

We always choose $\mu$ to be a self-similar measure on $K$ in this paper. To be more precise, we fix a weight vector $\{\mu_i\}_{i=1}^N$, and let $\mu$ be the unique probability measure supported on $K$ such that 
\[\mu(A)=\sum_{i=1}^N \mu_i\mu(F_i^{-1}A),\quad \forall A\subset K.\]
One can easily check that $\mu(F_wK)=\mu_w:=\mu_{w_1}\cdots \mu_{w_m}$, for each $w\in W_m$. 

For $f\in dom\mathcal{E}$, say $\Delta f=u$ if 
\[\mathcal{E}(f,\varphi)=-\int_K u\varphi d\mu \]
holds for any $\varphi\in dom_0\mathcal{E}$, with $dom_0\mathcal{E}=\{\varphi\in dom\mathcal{E}:\varphi|_{V_0}=0\}$. 

Write $L^2(K,\mu)=L^2(K)$ for short.

\begin{definition}\label{sobolev}
	For $k\in \mathbb{Z}_+$, define the Sobolev space $H^{2k}(K)$ as 
	\[H^{2k}(K)=\{f\in L^2(K): \Delta^j f\in L^2(K) \text{ for all } j\leq k\}\]
	with the norm $\|f\|_{H^{2k}(K)}$ of $f$ given by
	\[\|f\|^2_{H^{2k}(K)}=\sum_{j=0}^k \|\Delta^j f\|^2_{L^2(K)}\asymp \|f\|^2_{L^2(K)}+\|\Delta^k f\|^2_{L^2(K)}.\]
	For $0<\theta<1$, $k\in \mathbb{Z}_+$, define $H^{2k+2\theta}(K)$ to be
	\[H^{2(k+\theta)}(K)=[H^{2k}(K),H^{2k+2}(K)]_\theta\]
	the complex interpolation space.
\end{definition}

Analogously, by additionally requiring that each $\Delta^j f$ satisfies the Dirichlet boundary condition for $j<k$ in the above definition when {$k\in\mathbb{Z}_+$}, {we get a subspace, denoted by $H^{2k}_D(K)$, of $H^{2k}(K)$. The definition can be extended to any $\sigma\geq 0$ by using Bessel type potentials.} For $\sigma\geq 0$, we have $H^{\sigma}_D(K)=(id-\Delta_D)^{-\sigma/2}L^2(K)$, {with norm $\|(id-\Delta_D)^{\sigma/2}f\|_{L^2(K)}$}, where $\Delta_D$ is the Dirichlet Laplacian. {In particular, for $k\in \mathbb{Z}_+$ and $f\in H^{2k}_D(K)$, we have
\[\|f\|_{H^{2k}(K)}\asymp\|f\|_{H^{2k}_D(K)}\asymp \|\Delta^k f\|_{L^2(K)}.\]
Similarly,} we can define $H^{\sigma}_N(K)=(id-\Delta_N)^{-\sigma/2}L^2(K)$ with $\Delta_N$ being the Neumann Laplacian. See \cite{s1} by Strichartz for more details.

Throughout the paper, we always write $f\lesssim g$ if $f\leq Cg$ for some constant $C>0$, and write $f\asymp g$ if both $f\lesssim g$ and $g\lesssim f$.

\begin{comment}
Through section 2 to 6, we always assume\vspace{0.1cm}

(\textbf{A2}). {\color{red}\textit{There is $d_H>0$ such that $\mu_i=r_i^{d_H}$ for each $1\leq i\leq N$.}}\vspace{0.1cm}

We remark here that (\textbf{A2}) is not a necessary condition for most of our results, except those for low orders. See section 6 for further discussions for the case without the assumptions. Below we recall two useful results in our previous work, \cite{cq3}. 

\begin{proposition}\label{prop21}
	$H^1(K)=dom\mathcal{E}$.
\end{proposition}

\begin{proposition}\label{prop22}
	
\end{proposition}
\end{comment}

\section{Boundary behaviors of $H^\sigma(K)$}
In our previous work \cite{cq3}, we have made a complete comparison between $H^\sigma(K)$ and $H_D^\sigma(K)$ or $H_N^\sigma(K)$ for $\sigma\geq 0$. In this paper, we will make a more delicate description of the boundary behaviors of functions in $H^\sigma(K)$. For this purpose, a neat definition of tangents of functions at boundary points of $K$ is needed. 

It is natural to define tangents of functions as elements of multiharmonic functions. We will present our definition modified from that of Rogers and Strichartz \cite{distribution,s3}. 

\begin{definition}\label{def31}
For $k\geq 0$, let $\mathcal{H}_k=\{f\in H^{2k+2}(K):\Delta^{k+1}f=0\}$ be the space of $(k+1)$-multiharmonic functions on $K$. Let $\mathcal{H}_\bullet=\bigcup_{k=0}^\infty \mathcal{H}_k$.
\end{definition}  

Throughout this paper, we write $X=\oplus_{k=1}^n X_k$ for Banach spaces $X$ and $X_k,1\leq k\leq n$, if 

\noindent1. $X_k\subset X$ and $\|\cdot\|_{X_k}\asymp \|\cdot\|_X$, for each $1\leq k\leq n$;

\noindent2. For each $x\in X$, there is a unique representation $x=\sum_{k=1}^nx_k$, with $x_k\in X_k, 1\leq k\leq n$.

\begin{definition}\label{def32}  Fix $w\in W_*$.
	
(a). Define $A_w$ by $A_wf(x)=f(F_wx)$ for any function $f$ on $K$.

(b). Let {$\{\lambda_{l,w}\}_{l=0}^\infty$} be the set of nonzero eigenvalues of $A_w:\mathcal{H}_\bullet\to\mathcal{H}_\bullet$, which is ordered in decreasing order of absolute values, i.e. $1=|\lambda_{0,w}|\geq |\lambda_{1,w}|\geq |\lambda_{2,w}|\geq\cdots$. Let $E_{l,w}=\bigcup_{n=1}^\infty \ker(A_w-\lambda_{l,w})^n\subset \mathcal{H}_\bullet$ be the generalized eigenspace of $A_w$ corresponding to $\lambda_{l,w}$.

(c). Define $\tilde{E}_{l,w}=\oplus_{i=0}^{l''} E_{i,w}$ and  $\hat{E}_{l,w}=\oplus_{i=l'}^{l''} E_{i,w}$, where $${l'=\min\{i\geq 0:|\lambda_{i,w}|=|\lambda_{l,w}|\},l''=\max\{i\geq 0:|\lambda_{i,w}|=|\lambda_{l,w}|\}.}$$
\end{definition}

\noindent\textbf{Remark.} (a). It is well known that $\lim_{l\to\infty} |\lambda_{l,w}|=0$, and $E_{l,w}$ is of finite dimension for each $l$. In addition, we have $\lambda_{0,w}=1$ with $E_{0,w}=\hat{E}_{0,w}=\tilde{E}_{0,w}=constants$. 

(b). If $|\lambda_{l_1,w}|=|\lambda_{l_2,w}|$, then $\tilde{E}_{l_1,w}=\tilde{E}_{l_2,w}$ and $\hat{E}_{l_1,w}=\hat{E}_{l_2,w}$. In other words, the definition of $\tilde{E}_{l,w}$ and $\hat{E}_{l,w}$ only depends on the absolute value of $\lambda_{l,w}$.\vspace{0.1cm}

Now, we define the tangent of a function $f$ at a boundary vertex $p\in V_0$.

\begin{definition}\label{def33} Let $f\in C(K)$,  $\omega=\tau\dot{w}\in\mathcal{P}$ with $\dot{w}=ww\cdots$, and $l\geq 0$. {A multiharmonic function $h\in \tilde{E}_{l,w}$ is called a $l$-tangent of $f$ at $\pi(\omega)$ if 
	\[\|A_\tau f-h\|_{L^\infty(F_w^nK)}=o(\lambda_{l,w}^n),\]
and we denote $T_{l,\omega}f:=h$. In particular, $\|f-T_{l,\omega}f\|_{L^\infty(F_w^nK)}=o(\lambda_{l,w}^n)$ if $\tau=\emptyset$.}
\end{definition}

Intuitively, each $\omega\in \mathcal{P}$ represents a boundary point $p\in V_0$ and a ``direction'' that approaches $p$. We can view the collection of tangents $T_{l,\omega}, {\omega}\in\pi^{-1}(p)$ as the ``tangent'' at the boundary point $p\in V_0$. In the rest of this section, we will study the $l$-tangents for functions in $H^\sigma(K)$.

\subsection{Two lemmas}
Let $X$ be a Banach space, with norm $\|\cdot\|_X$, and $A:X\to X$ be a compact operator. Denote $\sigma(A,X)$  the spectrum of $A:X\to X$. We consider the following sequence spaces in this subsection. 

\begin{definition}\label{def34}
	(a). For $\alpha>0$, define $$l^2(X;\alpha)=\big\{\bm{s}=\{s_n\}_{n=0}^\infty:\{\alpha^{-n}\|s_n\|_X\}_{n=0}^\infty\in l^2\},$$
	with norm $\|\bm{s}\|_{l^2(X;\alpha)}=\big\|\alpha^{-n}\|s_n\|_X\big\|_{l^2}$.
	
	(b). For $\alpha>0$, define 
	\[l^2(X,A;\alpha)=\big\{\bm{s}=\{s_n\}_{n=0}^\infty:\{s_{n+1}-As_n\}_{n=0}^\infty\in l^2(X;\alpha)\big\},\]
	with norm $\|\bm{s}\|_{l^2(X,A; \alpha)}=\|s_{n+1}-As_n\|_{l^2(X;\alpha)}+\|s_0\|_X.$
	
	(c). For each $s\in X$,  define $\mathcal{S}_A(s)=\{A^ns\}_{n=0}^\infty$, with norm $\|\mathcal{S}_A(\cdot)\|_{\mathcal{S}_A(X)}=\|\cdot\|_{X}$.
\end{definition}

We let $\{\lambda_{l}\}_{l\geq 0}$ be the nonzero eigenvalues of $A$, which is ordered in decreasing order of absolute values. Let $E_l$ be the corresponding generalized eigenspaces. In addition, denote $\tilde{E}_l=\oplus_{i=0}^{l''} E_i$ and $\hat{E}_l=\oplus_{i=l'}^{l''}E_i$, where ${l'=\min\{i\geq 0:|\lambda_i|=|\lambda_l|\}, l''=\max\{i\geq 0:|\lambda_i|=|\lambda_l|\}.}$

\begin{lemma}\label{lemma35}
	 (a). For $\alpha>|\lambda_0|$ or $\sigma(A,X)=\{0\}$, we have $l^2(X,A;\alpha)=l^2(X;\alpha)$.
	
	(b). For $|\lambda_{l+1}|<\alpha<|\lambda_{l}|$ or $\alpha<|\lambda_l|=\min\{|\lambda_k|: \lambda_k\in\sigma(A,X)\}$, we have $$l^2(X,A;\alpha)=\mathcal{S}_A(\tilde{E}_l)\oplus l^2(X;\alpha).$$
\end{lemma}
\textit{Proof.} Let $\bm{s}\in l^2(X,A;\alpha)$, and denote $t_0=s_0$ and $t_n=s_{n}-As_{n-1}$ for $n\geq 1$. Clearly, $\bm{t}:=\{t_n\}_{n\geq 0}\in l^2(X;\alpha)$ with $\|\bm{t}\|_{l^2(X;\alpha)}\asymp\|\bm{s}\|_{l^2(X,A;\alpha)}$, and also $s_n=\sum_{m=0}^n A^{n-m}t_m$.

(a).  Using {the} Minkowski inequality, noticing that $\alpha$ is larger than the spectral radius of $A:X\to X$, we get
\[\begin{aligned}
\|\bm{s}\|_{l^2(X;\alpha)}&=\|\sum_{m=0}^n A^{n-m}t_m\|_{l^2(X;\alpha)}\leq\big\|\sum_{m=0}^n \alpha^{-n} \|A^{m}t_{n-m}\|_X\big\|_{l^2}\\
&\leq \big\|\sum_{m=0}^\infty 1_{n\geq m}\alpha^{m-n}\|\alpha^{-m} A^{m}  t_{n-m}\|_X\big\|_{l^2}\leq \sum_{m=0}^\infty \|\alpha^{-m}A^m\|_{X\to X} \cdot\|\bm{t}\|_{l^2(X;\alpha)}\\
&\lesssim \|\bm{s}\|_{l^2(X,A;\alpha)}.
\end{aligned}\]
 The other direction estimate is obvious.

(b). First, we assume $A$ is of finite rank and $X=\oplus_{i=0}^l E_i$, so that $A^{-1}$ is well defined. Clearly,  the following limit exists in $X$,
\[s_{\lim}=\lim_{n\to\infty} A^{-n}s_n=\sum_{m=0}^\infty A^{-m}t_m,\]
since $\alpha^{-1}$ is larger than the spectral radius of $A^{-1}$. Thus $$s_n-A^ns_{\lim}=-\sum_{m=n+1}^\infty A^{n-m}t_m.$$
Now we define $\bm{s}'=\{s_n-A^ns_{\lim}\}_{n\geq 0}$. Using Minkowski inequality, we get
\[\begin{aligned}
\|\bm{s}'\|_{l^2(X;\alpha)}&=\big\|\alpha^{-n}\|\sum_{m=n+1}^\infty A^{n-m}t_m\|_{X}\big\|_{l^2}\\
&\leq \big\|\sum_{m=1}^\infty \alpha^{-n-m}\| \alpha^{m}A^{-m}t_{n+m}\|_X\big\|_{l^2}\leq \sum_{m=1}^\infty \|\alpha^{m}A^{-m}\|_{X\to X}\cdot\|\bm{t}\|_{l^2(X;\alpha)}\\
&\lesssim \|\bm{s}\|_{l^2(X,A;\alpha)}.
\end{aligned}\]
This shows that $\bm{s}'\in l^2(X;\alpha)$ with the estimate of the norm. Clearly, the decomposition of $\bm{s}$ is unique, and both $\mathcal{S}_{A}(X)$ and $l^2(X;\alpha)$ are closed subspace of $X$. This proves (b) for the case that $A$ is of finite rank.

For general case, since $A$ admits a discrete spectrum, we can find a closed  subspace $\tilde{E}_l^\circ$ such that $\sigma(A,\tilde{E}^\circ_l)=\sigma(A,X)\setminus \{\lambda_i\}_{i=0}^l$, and $X=\tilde{E}_l\oplus \tilde{E}_l^\circ$. Then we see that
\[\begin{aligned}
l^2(X,A;\alpha)&=l^2(\tilde{E}_l,A;\alpha)\oplus l^2(\tilde{E}^{\circ}_l,A;\alpha)\\
&=\mathcal{S}_A(\tilde{E}_l)\oplus l^2(\tilde{E}_l;\alpha)\oplus l^2(\tilde{E}^{\circ}_l;\alpha)=\mathcal{S}_A(\tilde{E}_l)\oplus l^2(X;\alpha),
\end{aligned}\]
where we use (a) and the finite rank case that we have proved for (b).\hfill$\square$\vspace{0.2cm}

The case that $\alpha=|\lambda_l|$ for some $l\geq 0$ is a little more complicated. See the following lemma.
\begin{lemma}\label{lemma36} Let $\overline{l^2(X;\alpha)}$ be the closure of $l^2(X;\alpha)$ in $l^2(X,A;\alpha)$. Then $\overline{l^2(X;\alpha)}={l^2(X;\alpha)}$ if and only if $\alpha\notin\{|\lambda_{l}|\}_{l\geq 0}$. In addition,
	
(a). For $\alpha\geq |\lambda_0|$ or $\sigma(A,X)=\{0\}$, we have $l^2(X,A;\alpha)=\overline{l^2(X;\alpha)}$.

(b). For $|\lambda_{l+1}|\leq \alpha<|\lambda_{l}|$ or $\alpha<|\lambda_l|=\min\{|\lambda_k|: \lambda_k\in\sigma(A,X)\}$, we have $$l^2(X,A;\alpha)=\mathcal{S}_A(\tilde{E}_l)\oplus \overline{l^2(X;\alpha)}.$$
\end{lemma}
\textit{Proof.} 
(a). By Lemma \ref{lemma35} (a), we only need to prove the assertion for $\alpha=|\lambda_0|$. It suffices to show that $\mathcal{S}_A(\hat{E}_0)\subset \overline{l^2(X;\alpha)}$, since $\overline{\mathcal{S}_A(\hat{E}_0)+l^2(X;\alpha)}\supset \overline{l^2(X,A;\alpha-\varepsilon)}=l^2(X,A;\alpha)$ for some small $\varepsilon>0$ by Lemma \ref{lemma35} (b).

Let $\bm{s}=\mathcal{S}_A(s)$ for some $s\in E_0$. There is a $d\geq 0$ such that $(A-\lambda_0)^ds\neq 0$ and $(A-\lambda_0)^{d+1}s=0$. Write $s^{(k)}=(A-\lambda_0)^{k}s,0\leq k\leq d$. Fix $m_0,m_1,\cdots,m_d\in \mathbb{N}$ and take $M_k=\sum_{i=0}^{k} m_i$ (set $M_{-1}=0$). Then we can design a sequence $\bm{s}^{(m_0,m_1,\cdots,m_d)}=\{s_n^{(m_0,m_1,\cdots,m_d)}\}_{n\geq 0}$ in $l^2(X;\alpha)$ as follows. 
\[\lambda_0^{-n} s_n^{(m_0,m_1,\cdots,m_d)}=\begin{cases}
s,&\hspace{-3cm}\text{ if }n=0,\\
\lambda_0^{-n}As^{(m_0,m_1,\cdots,m_d)}_{n-1}-m_k^{-1}a_k s^{(k)},&\hspace{-3cm}\text{ if } M_{k-1}< n\leq M_k,\\
\text{ where } \lambda_0^{-M_{k-1}} s_{M_{k-1}}^{(m_0,m_1,\cdots,m_d)}=a_ks^{(k)}+b_{k+1}s^{(k+1)}+\cdots+b_ds^{(d)},\\
0, &\hspace{-3cm}\text{ if }n>M_d.
\end{cases}
\]
We can easily check that
\[\lim_{m_0\to\infty}\lim_{m_1\to\infty}\cdots \lim_{m_d\to\infty} \|\mathcal{S}_A(s)-\bm{s}^{(m_0,m_1,\cdots,m_d)}\|_{l^2(X,A;\alpha)}=0,\]
which gives that $\mathcal{S}_A(E_0)\subset \overline{l^2(X;\alpha)}$. Similarly, we have $\mathcal{S}_A(\hat{E}_0)\subset \overline{l^2(X;\alpha)}$.

The assertion (b) follows from a same argument as the proof of Lemma \ref{lemma35} (b).

As a consequence of Lemma \ref{lemma35}, we have $l^2(X;\alpha)=\overline{l^2(X;\alpha)}$ if  $\alpha\notin\{|\lambda_{l,w}|\}_{l\geq 0}$. On the other hand, if $\alpha=|\lambda_l|$ for some $l\geq 0$, we have $\mathcal{S}_A(s)\in \overline{l^2(X;\alpha)}\setminus l^2(X;\alpha)$ for $s\in E_l$.\hfill$\square$\vspace{0.1cm}

\noindent\textbf{Remark.} In the rest of this paper, without further clarification, we will always take $\overline{l^2(X;\alpha)}$ to be the closure of $l^2(X;\alpha)$ in $l^2(X,A;\alpha)$ as in Lemma \ref{lemma36}. This space will play an important role in Section 4 and 5.

\subsection{Construction of tangents: higher order case}
In this subsection, we construct the tangents for functions $f\in H^\sigma(K)$ at points in $V_0$ for $\sigma\geq 2$. For simplicity, we will fix a $\omega=\dot{w}\in\mathcal{P}$.

\noindent\textbf{Notation.} \textit{Let $X$ be a closed subspace of $L^2(K)$, and $w\in W_*$.}
	
	 \textit{(a).  Define $A_{w,X}=A_w\circ P_X$, where $P_X$ is the orthogonal projection from $L^2(K)$ onto $X$.}
	
	 \textit{(b). Define $\mathcal{A}_w(f)=\{A_w^nf\}_{n\geq 0}$ for each $f\in L^2(K)$.}

We start by studying functions in $H^{2k}(K),k\geq 1$. 
\begin{lemma}\label{lemma37}
Let $k\geq 1$. There exists $g_{w,k}\in C(K\times K)$ such that for any $f\in H^{2k}(K)$, we have 
\[A_wf(x)-A_{w,\mathcal{H}_{k-1}}f(x)=\int_{K} g_{w,k}(x,y)(-\Delta)^k f(y)d\mu(y).\]
\end{lemma}
\textit{Proof.} It is easy to see that $(A_w-A_{w,\mathcal{H}_{k-1}})|_{\mathcal{H}_{k-1}}=0$. So it suffices to prove the lemma for $f=G^k(-\Delta)^k f$, where $G$ is the Green's operator. 
We only need to take 
$$g_{w,k}(x,y)=\int_{K^{k-1}} (A_w-A_{w,\mathcal{H}_{k-1}})G_{y_1}(x)G(y_1,y_2)\cdots G(y_{k-1},y)d\mu({y_1})\cdots d\mu({y_{k-1}}),$$
where $G_y(x)=G(x,y)$ is the Green's function.\hfill$\square$

\begin{lemma}\label{lemma38}
	Let $f\in L^2(K)$ and $g\in L^\infty(K)$. Define $\bm{s}=\{\mu_w^{n/2}\int_{K} A^n_wf(x)g(x)d\mu(x)\}_{n= 0}^\infty$, then we have 
	\[\|\bm{s}\|_{l^2}\lesssim \|f\|_{L^2(K)}\|g\|_{L^\infty(K)}.\]
\end{lemma}
\textit{Proof.} Let $Z=K\setminus F_wK$. Then
 $\|f\|_{L^2(K)}=\big\|\|\mu_w^{n/2}A^n_wf\|_{L^2(Z)}\big\|_{l^2}$ by scaling, and 
 \[\begin{aligned}
 \big|\int_{K} f(x)g(x)d\mu(x)\big|&=\big|\sum_{m=0}^\infty\mu_w^{m}\int_{Z} A^m_wf(x)A^m_wg(x)d\mu(x)\big|\\
 &\leq \sum_{m=0}^\infty \mu_w^m\|A_w^mf\|_{L^2(Z)}\|g\|_{L^\infty(K)}.
 \end{aligned}\]
So using Minkowski inequality, we get
\[\begin{aligned}
\|\bm{s}\|_{l^2}&\leq \|g\|_{L^\infty(K)}\big\|\sum_{m=0}^\infty \mu_w^{n/2}\mu_w^m\|A_w^{n+m}f\|_{L^2(Z)}\big\|_{l^2}\\
& \leq \|g\|_{L^\infty(K)}\sum_{m=0}^\infty \mu_w^{m/2}\|f\|_{L^2(K)}.
\end{aligned}\]
Since $\mu_w<1$, we get the lemma. \hfill$\square$\vspace{0.2cm}

Using Lemma \ref{lemma37} and \ref{lemma38}, we get the following key observation.

\begin{proposition}\label{prop39}
	$\mathcal{A}_w\in \mathcal{L}\Big(H^\sigma(K),l^2\big(L^\infty(K),A_{w,\mathcal{H}_{k-1}};r_w^{\sigma/2}\mu_w^{(\sigma-1)/2}\big)\Big)$ for  $\sigma\geq 2$ and $k\geq \lceil \sigma/2\rceil$.
\end{proposition}
\textit{Proof.} First, we consider the $H^{2k}(K),k\geq 1$ case. We have the estimate 
\[\begin{aligned}
\|A_w^{n+1}f-A_{w,\mathcal{H}_{k-1}}A_w^nf\|_{L^\infty(K)}&\leq \big|\int_K \|g_{w,k}(\cdot,y)\|_{L^\infty(K)}(-\Delta)^k(A_w^nf)(y)d\mu(y)\big|\\
&=\big|r_w^{kn}\mu_w^{kn}\int_K \|g_{w,k}(\cdot,y)\|_{L^\infty(K)}A_w^n\big((-\Delta)^kf\big)(y)d\mu(y)\big|,
\end{aligned}\]
by using Lemma \ref{lemma37} and scaling.  Since $\|\Delta^k f\|_{L^2(K)}\leq \|f\|_{H^{2k}(K)}$, by using Lemma \ref{lemma38}, we have proved the assertion for $H^{2k}(K)$ cases.

For general case, we need to use the complex interpolation. For any $k'\geq k$, we see that
\[\begin{aligned}
l^2\big(L^\infty(K),A_{w,\mathcal{H}_{k-1}};r_w^{k}\mu_w^{k-1/2}\big)&=l^2\big(\mathcal{H}_{k-1},A_{w,\mathcal{H}_{k-1}};r_w^{k}\mu_w^{k-1/2}\big)+ l^2(L^\infty(K);r_w^{k}\mu_w^{k-1/2})\\
&=l^2\big(\mathcal{H}_{k-1},A_{w,\mathcal{H}_{k'}};r_w^{k}\mu_w^{k-1/2}\big)+ l^2(L^\infty(K);r_w^{k}\mu_w^{k-1/2})\\
&=l^2\big(L^\infty(K),A_{w,\mathcal{H}_{k'}};r_w^{k}\mu_w^{k-1/2}\big),
\end{aligned}\] 
where the first and last equalities are variants of Lemma \ref{lemma35}, using the fact that $\tilde{E}_{l,w}\subset \mathcal{H}_{k-1}$ for the largest $l$ with $|\lambda_{l,w}|\geq r_w^k\mu_w^{k-1/2}$.

The proposition then follows from the fact that
$$[l^2\big(X,A;\alpha),l^2\big(X,A;\beta)]_\theta=l^2\big(X,A;\alpha^{(1-\theta)}\beta^\theta),$$
and $[H^{2k}(K),H^{2k+2}(K)]_\theta=H^{2k+2\theta}(K)$ with $\theta\in[0,1]$. \hfill$\square$\vspace{0.2cm}

{Using Lemma \ref{lemma36}}, we get the existence of  tangents for functions in $H^\sigma(K)$ with higher orders. 
 
\begin{comment}
\begin{theorem}\label{thm1}
	Let $p=\pi(\tau\dot{w})\in V_0$ and $\lambda_{w,l}\vee \mu_w^{1/2}r_w<\mu_w^{(\sigma-1)/2}r_w^{\sigma/2}\leq \lambda_{w,l+1}$. We have $l-$tangent at $p$ for any $f\in H^\sigma(K)$. In addition we have 
    \[
    \begin{aligned}
    \sum_{n=0}^\infty \mu_w^{(1-\sigma)}r_w^{-\sigma}\|f-T_{l,p}f\|^2_{L^\infty(F_w^nK)}<\infty, 
    \end{aligned}
    \]
    if $\mu_w^{(\sigma-1)/2}r_w^{\sigma/2}\neq \lambda_{w,l+1}$.
\end{theorem}
\end{comment}

\subsection{Construction of tangents: lower order case}
The lower order case is a little more complicated. We still fix a $\omega=\dot{w}\in\mathcal{P}$. Intuitively, we would like to expect that when $r_w^{\sigma/2}\mu_w^{(\sigma-1)/2}<1=|\lambda_{0,w}|$, there exist some tangents of functions in $H^\sigma(K)$ at $p=\pi(\omega)$. However, to get this, we first need to guarantee that ${H}^\sigma(K)\subset L^\infty(K)$. Throughout this subsection and Section 4-6, we always assume the following assumption on  $\mu$.\vspace{0.1cm}

\noindent(\textbf{A2}) \textit{There exists $d_H>0$ such that $\mu_i=r_i^{d_H}$, $\forall 1\leq i\leq N$.}\vspace{0.1cm}

\noindent\textbf{Remark 1.} The assumption (\textbf{A2}) means that the self-similar measure $\mu$ is $d_H$-regular with respect to the effective resistance metric. See the authors' previous work \cite{cq3} for a discussion on the Besov characterizations of $H^\sigma(K)$ under (\textbf{A2}). We will use some results from \cite{cq3}.

\noindent\textbf{Remark 2.} Recall that the \textit{spectral dimension} of $K$ is $d_S=\frac{2d_H}{1+d_H}$. Then $\frac{d_S}{2}$ is the critical value of $\sigma$ such that $r_i^{\sigma/2}\mu_i^{(\sigma-1)/2}=1$ for all $1\leq i\leq N$. {In addition, we have $H^\sigma(K)\subset L^\infty(K)$ if and only if $\sigma>\frac{d_S}{2}$, see Theorem 6.1 in \cite{cq3}.} When $K$ is the unit interval, we have $\frac{d_S}{2}=\frac{1}{2}$, which is indeed the critical order in the Euclidean case. See \cite{Lion}. 

For $0<t\leq 1$, define $$\Lambda(t)=\{u\in W_*:r_u\leq t<r_{u^*}\},$$ where for $u=u_1u_2\cdots u_m$, $u^*=u_1u_2\cdots u_{m-1}$.
In particular, set $r=\min_{i=1}^N r_i$, and let $\Lambda_m=\Lambda(r^m)$ for $m\geq 0$. For any $u\in W_*$, define the average of $f$ on $F_uK$ by $Avg_u(f)=\mu_u^{-1}\int_{F_uK}fd\mu.$
In particular, $Avg_\emptyset(f)=\int_{K}fd\mu$.

For $m\geq 1$, define the space of \textit{$m$-Haar functions} 
\[\tilde{J}_m=\{\tilde{f}_m=\sum_{u\in \Lambda_m} c_u1_{F_uK}: c_u\in \mathbb{R}, Avg_{u'}(\tilde{f}_m)=0,\forall u'\in \Lambda_{m-1}\},\]
where $1_{E}$ is the \textit{characteristic function} of a set $E$. Let $J_0$ be the space of constant functions. Let $$\tilde{\Gamma}_\sigma(K)=\big\{f\in L^2(K):\{P_{\tilde{J}_n}f\}_{n=0}^\infty \in l^2(L^2(K); r^{\sigma(1+d_H)/2})\big\},$$
with norm $\|\cdot\|_{\tilde{\Gamma}_\sigma(K)}=\|\{P_{\tilde{J}_n}\cdot\}_{n=0}^\infty\|_{l^2(L^2(K);r^{\sigma(1+d_H)/2})}$. The following result comes from Theorem 3.9 and 4.11  in \cite{cq3}.

\begin{proposition}\label{thm310}
	 For $\frac{d_S}{2}<\sigma<1$, we have $H^\sigma(K)=\tilde{\Gamma}_\sigma(K)\cap C(K)$, with $\|f\|_{H^\sigma(K)}\asymp \|f\|_{\tilde{\Gamma}_\sigma(K)}$.
\end{proposition}

Using Proposition \ref{thm310}, we get the following estimate.
\begin{lemma}\label{lemma311}
	For each $f\in H^\sigma(K)$ with $\frac{d_S}{2}<\sigma<1$, let $f_n=\sum_{m=n}^\infty P_{\tilde{J}_m}f$. Then, we have $\{f_n\}_{n\geq 0}\in l^2(L^\infty(K);r^{\sigma/2}r^{(\sigma-1)d_H/2})$.
\end{lemma}

\textit{Proof.} Noticing that $\|P_{\tilde{J}_n}f\|_{L^\infty(K)}\lesssim r^{-nd_H/2}\|P_{\tilde{J}_n}f\|_{L^2(K)}$, by using Minkowski inequality, we have
\[\begin{aligned}
&\big\|r^{-n\sigma/2}r^{n(1-\sigma)d_H/2}\|f_n\|_{L^\infty(K)}\big\|_{l^2}\lesssim \big\|r^{-n\sigma/2}r^{n(1-\sigma)d_H/2}\sum_{m=0}^\infty r^{-(m+n)d_H/2}\|P_{\tilde{J}_{m+n}}f\|_{L^2(K)}\big\|_{l^2}\\
=&\big\|\sum_{m=0}^\infty r^{m\sigma/2}r^{m(\sigma-1)d_H/2}  r^{-(m+n)\sigma(1+d_H)/2}\|P_{\tilde{J}_{m+n}}f\|_{L^2(K)}\big\|_{l^2}\\
\lesssim &\sum_{m=0}^\infty r^{m\sigma/2}r^{m(\sigma-1)d_H/2}\|f\|_{H^\sigma(K)}.
\end{aligned}\]
This finishes the proof. \hfill$\square$\vspace{0.2cm}

Using Lemma \ref{lemma311}, we can easily get the following proposition.

\begin{proposition}\label{prop312}
	$\mathcal{A}_w\in \mathcal{L}\Big(H^{\sigma}(K),l^2\big(L^\infty(K),A_{w,\mathcal{H}_{k-1}};r_w^{\sigma/2}\mu_w^{(\sigma-1)/2}\big)\Big)$  for  $\sigma>\frac{d_S}{2}$ and $k\geq \lceil \sigma/2\rceil$.
\end{proposition} 
\textit{Proof}. By Proposition \ref{prop39}, it suffices to prove the case $\sigma<2$. For $\frac{d_S}{2}<\sigma<1$, it is easy to see that for $f\in H^{\sigma}(K)$ and $n\geq 0$, choosing $l$ such that $r^{l+1}< r_w^n\leq r^{l}$, we have
\[\|A^{n+1}_wf-A_wP_{\tilde{J}_0}A_w^nf\|_{L^\infty(K)}\leq \|A^n_wf-P_{\tilde{J}_0}A_w^nf\|_{L^\infty(K)}\lesssim \|\sum_{m=l+1}^\infty P_{\tilde{J}_m}f\|_{L^\infty(K)}.\]
It then follows from Lemma \ref{lemma311} that 
$$\{A_w^nf\}_{n\geq 0}\in l^2(L^\infty(K), A_wP_{\tilde{J}_0};r_w^{\sigma/2}\mu_w^{(\sigma-1)/2}).$$
Using the well-known fact that $E_{0,w}=\tilde{J}_0$, the fact that $\lambda_{1,w}=r_w<r_w^{\sigma/2}\mu_w^{(\sigma-1)/2}<1$ and Lemma \ref{lemma35} (b), we see that 
\[l^2(L^\infty(K), A_wP_{\tilde{J}_0};r_w^{\sigma/2}\mu_w^{(\sigma-1)/2})=l^2(L^\infty(K), A_{w,\mathcal{H}_0};r_w^{\sigma/2}\mu_w^{(\sigma-1)/2}).\]
This implies that
\[\mathcal{A}_w\in \mathcal{L}\Big(H^\sigma (K),l^2\big(L^\infty(K),A_{w,\mathcal{H}_0};r_w^{\sigma/2}\mu_w^{(\sigma-1)/2}\big)\Big).\] 
The rest of the proof follows from complex interpolation, and using Proposition \ref{prop39}.\hfill$\square$

\subsection{A theorem of boundary behavior}Now, we conclude our results in the following theorem. For convenience, in the remaining of this paper, we use the following notations. 

\begin{definition}\label{def313}
	For $\sigma>d_S/2$ and $\omega=\tau\dot{w}\in \mathcal{P}$, let $l_\omega(\sigma)$ be the unique integer such that 
	$$|\lambda_{l_\omega(\sigma)+1,w}|\leq r_w^{\sigma/2}\mu_w^{(\sigma-1)/2}<|\lambda_{l_\omega(\sigma),w}|.$$
	For convenience,  we write $T_\omega^{(\sigma)}=T_{l_\omega(\sigma),\omega}$ for short. 
	
	In particular, when $\sigma\leq d_S/2$, we let $l_\omega(\sigma)=-1$ and $T_{l_\omega(\sigma),\omega}=0$. Also, set $\mathcal{H}_{-1}=\{0\}$ for convenience.
\end{definition}

From Proposition \ref{prop312} and Lemma \ref{lemma35}, we have
\begin{theorem}\label{thm314}
Let $\omega=\tau\dot{w}\in \mathcal{P}$ and $\sigma>d_S/2$. Then $A_\tau T_\omega^{(\sigma)}\in \mathcal{L}(H^\sigma(K),\tilde{E}_{l_\omega(\sigma),w})$, and 
\[
\begin{aligned}
\sum_{n=0}^\infty r_w^{-\sigma n}\mu_w^{(1-\sigma)n}\|A_\tau f-T_\omega^{(\sigma)} f\|^2_{L^\infty(F_w^nK)}\lesssim \|f\|^2_{H^\sigma(K)}, 
\end{aligned}
\]
if $|\lambda_{l_\omega(\sigma)+1,w}|<r_w^{\sigma/2}\mu_w^{(\sigma-1)/2}<|\lambda_{l_\omega(\sigma),w}|$.
\end{theorem}

\noindent\textbf{Remark.} Theorem \ref{thm314} is still true for $\sigma\geq 2$ without the assumption (\textbf{A2}).

\begin{comment}
\begin{definition}
	For each closed subspace $X$ of $L^2(K)$, we define $P_X:L^2(K)\to X$ to be the orthogonal projection. In addition, define $T_{w,X}:L^2(K)\to \big(L^2(K)\big)^\mathbb{Z_+}$ as follows
	\[T_{w,X}(f)=\{P_XA_w^nf\}_{n\geq 1}.\]
\end{definition}

\begin{lemma}
	$T_{w,X}\in \mathcal{L}(L^2(K);\mu_w^{-1/2})$, where $\mathcal{L}(A,B)$ is the set of continuous linear maps from $A$ to $B$. 
\end{lemma}
\textit{Proof.}  It suffices to show that $T_{w,L^2(K)}$ is continuous. Fix $f\in L^2(K)$,

\begin{lemma}
	(a). For any $f\in H^{2k}(K)$, we have $\|P_{\mathcal{H}_{k-1}^\bot}f\|_{H^{2k}(K)}\lesssim \|\Delta^k f\|_{L^2(K)}$. 
	
	(b). Then, $T_{w,\mathcal{H}_{k-1}^\bot}\in \mathcal{L}(H^{2k}(K),l^2(H^{2k}(K);r_w^k\mu_{w}^{k-1/2}))$. 
\end{lemma}
\end{comment}

\section{The Sobolev spaces $H^\sigma_0(K)$}
We now proceed to study the Sobolev spaces $H^\sigma_0(K)$. In this section, we will make a full characterization of the relationship between $H^\sigma_0(K)$ and $H^\sigma(K)$ in terms of the boundary behavior of functions. We assume (\textbf{A1}) and (\textbf{A2}) throughout this section. \vspace{0.1cm}

\noindent\textbf{Remark.} (a). All the theorems in this paper are true without (\textbf{A1}). For convenience, and to raise the readability of the paper, we choose to assume (\textbf{A1}) throughout Section 4 to  6.

(b). Assumption (\textbf{A2}) is important for small $\sigma$. However, all of the results in the following are true for $\sigma\geq 2$ without  (\textbf{A2}). \vspace{0.1cm}

Recall that our domain is $\Omega=K\setminus V_0$ with boundary $V_0$. The space of smooth functions with compact support is defined as 
\[\mathcal{D}(\Omega)=\{f\in C_c(\Omega):\Delta^k f\in C_c(\Omega), \forall k\geq 0\},\]
where $C_c(\Omega)$ is the space of continuous functions with compact support in $\Omega$. See {\cite{distribution}} for basic properties of $\mathcal{D}(\Omega)$ in the fractal settings. Naturally, $\mathcal{D}(\Omega)$ is a subspace of $H^\sigma(K)$, $\forall\sigma\geq 0$.

\begin{definition}\label{def41}
	For $\sigma\geq 0$, define $H^\sigma_0(K)$ as the closure of $\mathcal{D}(\Omega)$ in $H^\sigma(K)$ with respect to the norm $\|\cdot\|_{H^\sigma(K)}$. 
\end{definition}

\begin{theorem}\label{thm42}
	For $\sigma\geq 0$, we have $H^\sigma_0(K)=\{f\in H^\sigma(K):T^{(\sigma)}_\omega(f)=0,\forall \omega\in \mathcal{P}\}$. In particular, $H^\sigma_0(K)=H^\sigma(K)$ if $\sigma\leq \frac{d_S}{2}$.
\end{theorem}

Both Theorem \ref{thm314} and \ref{thm42} have elegant analogues for domains $\Omega\in\mathbb{R}^n$ with good boundary.  Related development can be found in \cite{Lion} (Chapter 1, Section 9 and 10). 

In the rest of this section, we will focus on the proof of Theorem \ref{thm42}. Since $T^{(\sigma)}_\omega$ is continuous by Theorem \ref{thm314} and $T^{(\sigma)}_\omega|_{\mathcal{D}(\Omega)}=0$, we can easily show one direction of the theorem holds, i.e.,
\[H^\sigma_0(K)\subset\{f\in H^\sigma(K):T^{(\sigma)}_\omega(f)=0,\forall \omega\in \mathcal{P}\}.\]
In the next two subsections, we will provide the proof for the other direction. 

\subsection{Pretangents}
The tangents $T^{(\sigma)}_\omega$ for different Sobolev spaces $H^\sigma(K)$ are taken with different orders. We hope to get more information. 

\begin{definition}\label{def43}
Let $X$ be a Banach space with an operator $A:X\to X$.  Raise $A$ to be $A:X^{\mathbb{Z}_+}\to X^{\mathbb{Z}_+}$ by  $$A\{s_n\}_{n\geq 0}=\{As_n\}_{n\geq 0}.$$ 
\end{definition} 

The following lemma is an easy observation from the last section.
\begin{lemma}\label{lemma44}
For $w\in W_*$, $\sigma\geq 0$ and $k\geq \lceil \sigma/2\rceil$, $P_{\mathcal{H}_{k-1}}\mathcal{A}_w\in \mathcal{L}\big(H^\sigma(K),l^2(\mathcal{H}_{k-1}, A_w;r_w^{\sigma/2}\mu_w^{(\sigma-1)/2})\big)$.
\end{lemma}
\textit{Proof.} For $\sigma=0$, it is easy to see the assertion by using Lemma \ref{lemma38}. In fact, $\mathcal{H}_{k-1}(K)$ is a finite dimensional subspace of $C(K)$, so there is $g_k\in C(K\times K)$ such that $P_{\mathcal{H}_{k-1}}f(x)=\int_K g_k(x,y)f(y)d\mu(y)$. Then a similar argument as in the proof of Proposition \ref{prop39} works.

For $\sigma\geq 2$ cases, the lemma is a consequence of Proposition \ref{prop39}, noticing that 
$$P_{\mathcal{H}_{k-1}}A_w^{n+1}f-A_wP_{\mathcal{H}_{k-1}}A_w^nf=P_{\mathcal{H}_{k-1}}\big(A_w^{n+1}f-A_{w,\mathcal{H}_{k-1}}A_w^{n}f\big).$$

For $0<\sigma<2$, the assertion can be proved by using complex interpolation. \hfill$\square$

\begin{theorem}\label{thm45}
For $\omega=\dot{w}\in\mathcal{P}$, $\sigma\geq 0$ and $k\in \mathbb{Z}^+$, there is a recovering map 
$$\mathcal{R}_{k,w}\in \mathcal{L}\big(l^2(\mathcal{H}_{k-1}, A_w;r_w^{\sigma/2}\mu_w^{(\sigma-1)/2}),H^\sigma(K)\big)$$
such that $P_{\mathcal{H}_{k-1}}\mathcal{A}_w\mathcal{R}_{k,w}=id.$ In addition, $\mathcal{R}_{k,w}(\cdot)$ vanishes in a neighbourhood of $V_0\setminus \{\pi(\omega)\}$.
\end{theorem}

\textit{Proof.} First we assume that $F_wK$ is bounded away from $V_0\setminus \{\pi(\omega)\}$. Then, for each $h\in \mathcal{H}_{k-1}$, obviously there is $f\in dom\Delta^\infty$ such that $A_wf=h, P_{\mathcal{H}_{k-1}}f=0$ and $f$ vanishes in a neighbourhood of $V_0\setminus \{\pi(\omega)\}$. By a standard argument,  there is a {linear} map $R_{k,w}:\mathcal{H}_{k-1}\to dom\Delta^\infty$ such that $A_wR_{k,w}(h)=h$, $P_{\mathcal{H}_{k-1}}R_{k,w}(h)=0$ and $R_{k,w}(h)$ vanishes in a neighbourhood of $V_0\setminus \{\pi(\omega)\}$.

By a same reason, there is a map $R'_{k,w}:\mathcal{H}_{k-1}\to dom\Delta^\infty$ such that  $A_wR'_{k,w}(h)=A_wh$, $P_{\mathcal{H}_{k-1}}R'_{k,w}(h)=h$ and $R'_{k,w}(h)$ vanishes in a neighbourhood of $V_0\setminus \{\pi(\omega)\}$. 

Now for any $\bm{h}=\{h_n\}_{n\geq 0}\in l^2(\mathcal{H}_{k-1}, A_w;r_w^{\sigma/2}\mu_w^{(\sigma-1)/2})$, we define 
\[\mathcal{R}_{k,w}(\bm{h})=R'_{k,w}(h_0)+\sum_{n=1}^\infty R_{k,w}(h_n-A_wh_{n-1})\circ F^{-n+1}_w.\]
We need to show that $\mathcal{R}_{k,w}$ is well defined and is in $\mathcal{L}\big(l^2(\mathcal{H}_{k-1}, A_w;r_w^{\sigma/2}\mu_w^{(\sigma-1)/2}),H^\sigma(K)\big)$. 

First, we consider $\sigma=2k'<2k$. Write  $f_n=R_{k,w}(h_n-A_wh_{n-1})\circ F_{w}^{-n+1},n\geq 1$ and $f_0=R'_{k,w}(h_0)$ for short. Then we see that 
$$\begin{aligned}
\|\{\Delta^{k'} f_n\}_{n\geq 0}\|_{l^2(L^\infty(K);\mu_w^{-1/2})}&\lesssim \|h_0\|_{\mathcal{H}_{k-1}}+\|\{h_n-A_wh_{n-1}\}_{n\geq 0}\|_{l^2(\mathcal{H}_{k-1};r_w^{k'}\mu_w^{k'-1/2})}\\&=\|\bm{h}\|_{l^2(\mathcal{H}_{k-1},A_w;r_w^{k'}\mu_w^{k'-1/2})}.
\end{aligned}$$
Write $Z=K\setminus F_wK$. Then we have 
\[
\begin{aligned}
&\big\|\sum_{m=0}^\infty |\Delta^{k'}f_m|\big\|_{L^2(K)}=\big\|\mu^{n/2}_w\|A^n_w\sum_{m=0}^{n+1} |\Delta^{k'}f_m|\|_{L^2(Z)}\big\|_{l^2}\lesssim \big\|\mu^{n/2}_w\sum_{m=0}^{{ n+1}}\|\Delta^{k'}f_m\|_{L^\infty(K)}\big\|_{l^2}\\
=&\big\|{\sum_{m=-1}^\infty 1_{n\geq m}}\mu^{m/2}_w\mu^{(n-m)/2}_w\|\Delta^{k'}f_{n-m}\|_{L^\infty(K)}\big\|_{l^2}\lesssim \|\{\Delta^{k'}f_n\}_{n\geq 0}\|_{l^2(L^\infty(K);\mu_w^{-1/2})}.
\end{aligned}
\]
{Thus, $\big\|\sum_{m=0}^\infty|\Delta^{k'}f_m|\big\|_{L^2(K)}\lesssim \|\bm{h}\|_{l^2(\mathcal{H}_{k-1},A_w;r_w^{k'}\mu_w^{k'-1/2})}$. By a same argument, we have 
\begin{equation}\label{eqn41}\big\|\sum_{m=0}^\infty|f_m|\big\|_{L^2(K)}\lesssim \|\bm{h}\|_{l^2(\mathcal{H}_{k-1},A_w;\mu_w^{-1/2})}\leq \|\bm{h}\|_{l^2(\mathcal{H}_{k-1},A_w;r_w^{k'}\mu_w^{k'-1/2})}.
\end{equation}
The above estimates show that $\mathcal{R}_{k,w}(\bm{h})=\sum_{m=0}^\infty f_m$ is well defined in $H^{2k'}(K)$, and   $\mathcal{R}_{k,w}\in\mathcal{L}\big(l^2(\mathcal{H}_{k-1}, A_w;r_w^{k'}\mu_w^{k'-1/2}),H^{2k'}(K)\big)$. }

Next, we consider $\sigma=2k'\geq 2k$. Then, we can see that for any $n\geq 1$, $\Delta^{k'} f_n$ is supported in $F_w^{n-1}Z$  with $\|\Delta^{k'} f_n\|_{L^2(K)}\lesssim r_w^{-k'n}\mu_w^{-(k'-1/2)n}\|h_n-A_wh_{n-1}\|_{\mathcal{H}_{k-1}}$. This shows that 
\begin{equation}\label{eqn42}
\|\Delta^{k'}\mathcal{R}_{k,w}(\bm{h})\|_{L^2(K)}\lesssim \|\bm{h}\|_{l^2(\mathcal{H}_{k-1}, A_w;r_w^{k'}\mu_w^{k'-1/2})}.
\end{equation}
Combining estimates (\ref{eqn41}) and (\ref{eqn42}), we still see that $\mathcal{R}_{k,w}\in\mathcal{L}\big(l^2(\mathcal{H}_{k-1}, A_w;r_w^{k'}\mu_w^{k'-1/2}),H^{2k'}(K)\big)$. 

Using complex interpolation, we see that for any $\sigma\geq 0$, we have $$\mathcal{R}_{k,w}\in\mathcal{L}\big(l^2(\mathcal{H}_{k-1}, A_w;r_w^{\sigma/2}\mu_w^{(\sigma-1)/2}),H^\sigma(K)\big).$$
Lastly, it is easy to check that $P_{\mathcal{H}_{k-1}}\mathcal{A}_w\mathcal{R}_{k,w}=id$ from the definition.

For the case $F_wK\cap \big(V_0\setminus \{\pi(\omega)\}\big)\neq \emptyset$, we only need slightly modify the definition of $R_{k,w}$ and $R'_{k,w}$.\hfill$\square$\vspace{0.15cm}

To take care of all the boundary points at the same time, including those with addresses $\tau\dot{w}$ and $\tau\neq \emptyset$, we need to localize the recovering map. See the definition below.

\begin{definition}\label{def46}  Fix $l\geq 0$ such that for any distinct $\omega=\tau\dot{w}$ and $\omega'=\tau'\dot{w}'$ in $\mathcal{P}$, we have $F_\tau F_w^lK\cap F_{\tau'} F_{w'}^lK=\emptyset$.  Let $\omega=\tau \dot{w}\in \mathcal{P}$, $f$ be a function on $K$ and $\bm{h}\in \mathcal{H}_{k-1}^{\mathbb{Z}^+}$.
	
(a). Define $\mathcal{T}_{\mathcal{H}_{k-1},\omega}f=P_{\mathcal{H}_{k-1}}\mathcal{A}_w(A_w^{l}A_\tau f)$. Call $\mathcal{T}_{\mathcal{H}_{k-1},\omega}f$ the $k$-pretangent of $f$ at $\pi(\omega)$.

(b). Write $\mathcal{R}_{\mathcal{H}_{k-1},\omega}(\bm{h})=\mathcal{R}_{k,w}(\bm{h})\circ F_w^{-l}F_\tau^{-1}$.

(c). Define $ker_\sigma \mathcal{T}_{\mathcal{H}_{k-1}}=\{f\in H^\sigma(K):\mathcal{T}_{\mathcal{H}_{k-1},\omega}f=0, \forall \omega\in \mathcal{P}\}$, with {induced norm} $\|\cdot\|_{H^\sigma(K)}$. 
\end{definition}

\begin{theorem}\label{thm47}
Let $\sigma\geq 0$, $k\geq \lceil \sigma/2\rceil$.

(a). For each $\omega=\tau\dot{w}\in \mathcal{P}$ and $\sigma\geq 0$, we have \[\|\mathcal{R}_{\mathcal{H}_{k-1},\omega}(\bm{h})\|_{H^\sigma(K)}\asymp\|\bm{h}\|_{l^2(\mathcal{H}_{k-1},A_w;r_w^{\sigma/2}\mu_w^{(\sigma-1)/2})}.\]
	
(b).  $H^\sigma(K)=\ker_\sigma\mathcal{T}_{\mathcal{H}_{k-1}}\oplus\big(\oplus_{\omega\in \mathcal{P}}\mathcal{R}_{\mathcal{H}_{k-1},\omega}\big(l^2(\mathcal{H}_{k-1},A_w;r_w^{\sigma/2}\mu_w^{(\sigma-1)/2})\big)\big)$, $\forall \sigma\geq 0$.

(c). Let $\sigma'>\sigma\geq 0$, then $ker_{\sigma'} \mathcal{T}_{\mathcal{H}_{k-1}}$ is a dense subspace of $ker_{\sigma} \mathcal{T}_{\mathcal{H}_{k-1}}$.
\end{theorem}
\textit{Proof.} {(a) is an easy consequence of Lemma \ref{lemma44}, Theorem \ref{thm45}.}

(b). Given any function $f\in H^\sigma(K)$, we can easily see that 
\[f-\sum_{\omega\in\mathcal{P}} \mathcal{R}_{\mathcal{H}_{k-1},\omega}\circ\mathcal{T}_{\mathcal{H}_{k-1},\omega}(f)\in\ker_{\sigma}\mathcal{T}_{\mathcal{H}_{k-1}}.\]
The decomposition of $f$ is obviously unique.

(c). For a subspace $X\subset {H}^\sigma(K)$, we write $\overline{X}$ for the closure of $X$ in ${H}^\sigma(K)$ for short. For convenience, we also write $X_{\omega,\sigma}=\mathcal{R}_{\mathcal{H}_{k-1},\omega}\big(l^2(\mathcal{H}_{k-1},A_w;r_w^{\sigma/2}\mu_w^{(\sigma-1)/2})$ for short.

It is obvious that $\ker_{\sigma'} \mathcal{T}_{\mathcal{H}_{k-1}}\subset \ker_{\sigma} \mathcal{T}_{\mathcal{H}_{k-1}}$, and $X_{\omega,\sigma'}\subset X_{\omega,\sigma}$, which leads to
\[\begin{aligned}
\overline{H^{\sigma'}(K)}=\overline{\ker_{\sigma'} \mathcal{T}_{\mathcal{H}_{k-1}}}\oplus(\oplus_{\omega\in\mathcal{P}}\overline{X_{\omega,{\sigma'}}}).
\end{aligned}\]
However, we know that $\overline{H^{\sigma'}(K)}=H^{\sigma}(K)$ by a standard argument. This implies that $\overline{\ker_{\sigma'} \mathcal{T}_{\mathcal{H}_{k-1}}}=\ker_{\sigma} \mathcal{T}_{\mathcal{H}_{k-1}}$.
\hfill$\square$

\subsection{Proof of Theorem \ref{thm42}} Now we prove Theorem \ref{thm42}. The  \textit{smooth bump functions} developed by Rogers, Strichartz and Teplyaev in \cite{smoothbump} will play a key role in the proof. In particular, we will use the following easy consequence (See Theorem 4.3 and estimate (4.7) in \cite{smoothbump}).

\begin{proposition}\label{thm48}
Let $k\geq 1, p\in V_0$ and $f\in dom\Delta^\infty$. There is a function $g\in dom\Delta^\infty$ such that 
	\[\begin{cases}
	\|g\|_{H^{2k}(K)}\lesssim \|f\|_{H^{2k}(K)},\\
	\Delta^j g(q)=\Delta^j f(q), \partial_n \Delta^j g(q)=\partial_n \Delta^j f(q), \forall {j\geq 0}, \forall q\in V_0\setminus \{p\},
	\end{cases}\]
and the support of $g$ is away from $p$.	
\end{proposition}

Using Proposition \ref{thm48}, we have the following lemma.
\begin{lemma}\label{lemma49}
Let $k\geq 1$ and $f\in H^{2k}(K)$. If  $ \forall\omega=\tau\dot{w}\in\mathcal{P}$, $\mathcal{T}_{\mathcal{H}_{k-1},\omega}f\in l^2(\mathcal{H}_{k-1};r_w^k\mu_w^{k-1/2})$, then $f\in H^{2k}_0(K)$. 
\end{lemma}
\textit{Proof.} For each $\omega=\tau\dot{w}\in \mathcal{P}$, by an easy estimate, we can see that 
\[\|A^{n}_wA_\tau f-P_{\mathcal{H}_{k-1}}A^{n}_wA_\tau f\|_{H^{2k}(K)}\lesssim r_w^{kn}\mu_w^{kn}\|A^n_wA_\tau\Delta^k f\|_{L^2(K)}=o(r_w^{kn}\mu_w^{(k-1/2)n}).\] 
In addition, it is obvious from the assumption that
\[\|P_{\mathcal{H}_{k-1}}A^n_wA_\tau f\|_{H^{2k}(K)}\asymp \|P_{\mathcal{H}_{k-1}}A^n_wA_\tau f\|_{\mathcal{H}_{k-1}(K)}=o(r_w^{kn}\mu_w^{(k-1/2)n}).\]
Thus we have $\|A_w^nA_\tau f\|_{H^{2k}(K)}=o(r_w^{kn}\mu_w^{(k-1/2)n}).$ 

Now, we construct $g\in \mathcal{D}(\Omega)$ that well approximates $f$ in $H^{2k}(K)$. For any $\varepsilon>0$, we can do the following.

1. Choose a large $n$ such that $\|A_w^nA_\tau f\|_{H^{2k}(K)}\leq \varepsilon r_w^{kn}\mu_w^{(k-1/2)n}r_\tau^{k}\mu_\tau^{k-1/2}, \forall \omega\in \mathcal{P}$.

2. Let $P_t$ be the heat kernel generated by the Dirichlet Laplacian $\Delta_D$. Choose $t$ small enough so that $\|f-P_tf\|_{H^{2k}(K)}\leq \varepsilon$ and $\|A_w^nA_\tau (f-P_tf)\|_{H^{2k}(K)}\leq\varepsilon r_w^{kn}\mu_w^{(k-1/2)n}r_\tau^{k}\mu_\tau^{k-1/2}$. 

3. By Proposition \ref{thm48}, for each $\omega$, we can find a $g_\omega$ supported in $F_w^nF_\tau\setminus \pi(\omega)$ such that 
\[\begin{aligned}
\|\Delta^k g_\omega\|_{L^2(F_w^nF_\tau K)}&\leq r_w^{-kn}\mu_w^{-(k-1/2)n}r_\tau^{-k}\mu_\tau^{-k+1/2}\|A^n_wA_\tau g_\omega\|_{H^{2k}(K)}\\
&\leq Cr_w^{-kn}\mu_w^{-(k-1/2)n}r_\tau^{-k}\mu_\tau^{-k+1/2}\|A_w^nA_\tau P_t f\|_{H^{2k}(K)}\leq 2C\varepsilon,
\end{aligned}\] 
and $\Delta^j A_w^nA_\tau g_\omega(q)=\Delta^j A_w^nA_\tau P_tf(q), \partial_n \Delta^j A_w^nA_\tau g_\omega(q)=\partial_n\Delta^j A_w^nA_\tau P_tf(q),\forall {j\geq 0}, q\in V_0\setminus \{\pi(\omega)\}$.
Replace {$P_tf|_{F_w^nF_\tau(K)}$} with $g_\omega$ for each $\omega$, and name the induced function $g$. Clearly, $g\in \mathcal{D}(\Omega)$. One can then check that 
\[\|f-g\|_{H^{2k}(K)}\leq C'\|\Delta^kf-\Delta^k g\|_{L^2(K)}\leq C'(1+2\#\mathcal{P}+2C\#\mathcal{P})\varepsilon,\]
{noticing that both $f,g$ are in $H^{2k}_D(K)$. The lemma is proved by choosing $\varepsilon$ arbitrarily.}\hfill$\square$\vspace{0.2cm}

\textit{Proof of Theorem \ref{thm42}.} It suffices to show $\{f\in H^\sigma(K):T^{(\sigma)}_\omega(f)=0,\forall \omega\in \mathcal{P}\}\subset H^\sigma_0(K)$. Choose $k$ large enough so that $\sigma\leq 2k$. By Lemma \ref{lemma49}, we have 
\[\ker_{2k}\mathcal{T}_{\mathcal{H}_{k-1}}\oplus\Big(\oplus_{\omega\in \mathcal{P}}\mathcal{R}_{\mathcal{H}_{k-1},\omega}\big(l^2(\mathcal{H}_{k-1};r_w^{k}\mu_w^{k-1/2})\big)\Big)\subset H^{2k}_0(K)\subset H^{\sigma}_0(K).\]
As a consequence, by using Theorem \ref{thm47} (a) and (c), we get
\[\ker_\sigma\mathcal{T}_{\mathcal{H}_{k-1}}\oplus\Big(\oplus_{\omega\in \mathcal{P}}\mathcal{R}_{\mathcal{H}_{k-1},\omega}\big(\overline{l^2(\mathcal{H}_{k-1};r_w^{\sigma/2}\mu_w^{(\sigma-1)/2})}\big)\Big)\subset H^\sigma_0(K).\]
On the other hand, by Theorem \ref{thm45} and \ref{thm47} (b), and using Lemma \ref{lemma36}, we can see that
\[\{f\in H^\sigma(K):T^{(\sigma)}_\omega(f)=0,\forall \omega\in \mathcal{P}\}=\ker_\sigma\mathcal{T}_{\mathcal{H}_{k-1}}\oplus\Big(\oplus_{\omega\in\mathcal{P}}\mathcal{R}_{\mathcal{H}_{k-1},\omega}\big(\overline{l^2(\mathcal{H}_{k-1};r_w^{\sigma/2}\mu_w^{(\sigma-1)/2})}\big)\Big).\]
The assertion follows immediately.\hfill$\square$\vspace{0.1cm}

As a consequence of Theorem \ref{thm42}, we have the following characterization of $H^\sigma_0(K)$. 

\begin{theorem}\label{thm410}
	Let $k\geq 1$ be an integer, and let $0\leq\sigma\leq 2k$. Then we have 
	\begin{equation}\label{eqnthm410}
	H^\sigma_0(K)=\ker_\sigma\mathcal{T}_{\mathcal{H}_{k-1}}\oplus\Big(\oplus_{\omega\in \mathcal{P}}\mathcal{R}_{\mathcal{H}_{k-1},\omega}\big(\overline{l^2(\mathcal{H}_{k-1};r_w^{\sigma/2}\mu_w^{(\sigma-1)/2})}\big)\Big).
	\end{equation}
\end{theorem}

\section{Interpolation of $H^\sigma(K)$: $\sigma\geq 0$}
Now we are ready to turn to the second topic, the interpolation theorems. {In this section, we prove some interpolation theorems for Sobolev spaces with non-negative orders, and we will combine these results in a final theorem on Sobolev spaces with real orders in Section 6.} 

\begin{lemma}\label{lemma51}
Let $(Z_0,Z_1)$ be an interpolation couple, which means $Z_1$ and $Z_2$ are continuously embedded in a same Hausdoff topological vector space. Let $Z_0=X_0\oplus Y_0$ and $Z_1=X_1\oplus Y_1$, and assume that $(X_0+X_1)\cap (Y_0+Y_1)=\{0\}$. Then we have 

(a). $[Z_0,Z_1]_{\theta}=[X_0,X_1]_\theta\oplus [Y_0,Y_1]_\theta$;

(b). Assume that $[Z_0,Z_1]_{\theta}=\tilde{X}\oplus\tilde{Y}$ and $\tilde{X}\subset X_0+X_1, \tilde{Y}\subset Y_0+Y_1$, then $\tilde{X}=[X_0,X_1]_{\theta}$ and $\tilde{Y}=[Y_0,Y_1]_\theta$, with equivalent norms.
\end{lemma}
\textit{Proof.} (a). {Since $Z_0+Z_1=(X_0+X_1)+(Y_0+Y_1)$ with $(X_0+X_1)\cap (Y_0+Y_1)=\{0\}$, we can define the natural projection $P:Z_0+Z_1\to X_0+X_1$ such that $(1-P):Z_0+Z_1\to Y_0+Y_1$. It is easy to check that $P\in \mathcal{L}([Z_0,Z_1]_\theta,[X_0,X_1]_\theta)$ and $1-P\in \mathcal{L}([Z_0,Z_1]_\theta,[Y_0,Y_1]_\theta)$ by using complex interpolation. So $[Z_0,Z_1]_\theta=[X_0,X_1]_\theta+[Y_0,Y_1]_\theta$ with $[X_0,X_1]_\theta\cap[Y_0,Y_1]_\theta=\{0\}$. 

It remains to check that $\|x\|_{[Z_0,Z_1]_\theta}\asymp \|x\|_{[X_0,X_1]_\theta}$ for any $x\in [X_0,X_1]_\theta$ and $\|y\|_{[Z_0,Z_1]_\theta}\asymp \|y\|_{[Y_0,Y_1]_\theta}$ for any $y\in [Y_0,Y_1]_\theta$. Let $i_X$ be the embedding map from $X_0+X_1\to Z_0+Z_1$. Then, by using complex interpolation, one can see that $i_X\in \mathcal{L}([X_0,X_1]_\theta,[Z_0,Z_1]_\theta)$. As a consequence, $\|x\|_{[X_0,X_1]_\theta}=\|Px\|_{[X_0,X_1]_\theta}\lesssim \|x\|_{[Z_0,Z_1]_\theta}=\|i_Xx\|_{[Z_0,Z_1]_\theta}\lesssim \|x\|_{[X_0,X_1]_\theta}$, where we also use the fact that $P\in \mathcal{L}([Z_0,Z_1]_\theta,[X_0,X_1]_\theta)$. The proof for $[Y_0,Y_1]_\theta$ is the same. }

(b). Clearly $[X_0,X_1]_\theta=P([Z_0,Z_1]_\theta)=\tilde{X}$, and $[Y_0,Y_1]_\theta=(1-P)([Z_0,Z_1]_\theta)=\tilde{Y}$. The estimate of the norms is obvious.\hfill$\square$

\begin{lemma}\label{lemma52}
	$\ker_\sigma \mathcal{T}_{\mathcal{H}_{k-1}},\sigma\geq 0$ (defined in Definition \ref{def46} (c)) is stable under complex interpolation, i.e., for $2k\geq\sigma>\sigma'\geq 0$, it holds that
	$$[\ker_\sigma \mathcal{T}_{\mathcal{H}_{k-1}},\ker_{\sigma'} \mathcal{T}_{\mathcal{H}_{k-1}}]_{\theta}=\ker_{(1-\theta)\sigma+\theta\sigma'} \mathcal{T}_{\mathcal{H}_{k-1}}, \forall \theta\in[0,1].$$
\end{lemma}
\textit{Proof.} We need to use the fact that $H^\sigma_D(K)$ is stable under complex interpolation, see {\cite{s1}}. In other words,  $H^{(1-\theta)\sigma+\theta\sigma'}_D(K)=[H^\sigma_D(K),H^{\sigma'}_D(K)]_\theta, \forall\theta\in[0,1]$. It is easy to see that $H^\sigma_D(K)=\ker_\sigma \mathcal{T}_{\mathcal{H}_{k-1}}\oplus X_\sigma$, where we write $X_\sigma=\oplus_{\omega\in\mathcal{P}} \mathcal{R}_{\mathcal{H}_{k-1},\omega}\mathcal{T}_{\mathcal{H}_{k-1},\omega}(H^\sigma_D(K))$ for convenience. Also, one can see
$\big(\ker_\sigma \mathcal{T}_{\mathcal{H}_{k-1}}+\ker_{\sigma'} \mathcal{T}_{\mathcal{H}_{k-1}}\big)\cap (X_\sigma+X_{\sigma'})=\ker_{\sigma'} \mathcal{T}_{\mathcal{H}_{k-1}}\cap X_{\sigma'}=\{0\}$, $\ker_{(1-\theta)\sigma+\theta\sigma'}\mathcal{T}_{\mathcal{H}_{k-1}}\subset \ker_{\sigma'} \mathcal{T}_{\mathcal{H}_{k-1}}$ {and $X_{(1-\theta)\sigma+\theta\sigma'}\subset X_{\sigma'}$}. The lemma follows from Lemma \ref{lemma51} (b).\hfill$\square$ 

\begin{lemma}\label{lemma53}
	Let $X$ be a Banach space and $A$ be a compact operator in $\mathcal{L}(X,X)$. Denote $\overline{l^2(X;\alpha)}$ the closure of $l^2(X;\alpha)$ in $l^2(X,A;\alpha)$. Then, we have 
	\[[\overline{l^2(X;\alpha)},\overline{l^2(X;\beta)}]_\theta=l^2(X;\alpha^{(1-\theta)}\beta^\theta),\]
	for any $\infty>\alpha>\beta>0$ and $0<\theta<1$.
\end{lemma}

\textit{Proof.} Let's first make some observations of the special cases. \vspace{0.1cm}

\textit{Claim 1:  Let $\beta=|\lambda_l|$, then $[\overline{l^2(E_l;\alpha)},\overline{l^2(E_l;\beta)}]_\theta=l^2(E_l;\alpha^{(1-\theta)}\beta^\theta)$. }

In this case, by using Lemma \ref{lemma35} (a) and Lemma \ref{lemma36}, we see that 
\[\begin{aligned}
{[}\overline{l^2(E_l;\alpha)},\overline{l^2(E_l;\beta)}{]}_\theta&=[l^2(E_l;\alpha),l^2(E_l,A;\beta)]_\theta=[l^2(E_l,A;\alpha),l^2(E_l,A;\beta)]_{\theta}\\
&=l^2(E_l,A;\alpha^{(1-\theta)}\beta^\theta)=l^2(E_l;\alpha^{(1-\theta)}\beta^\theta),
\end{aligned}\]
where the last equality holds since $\alpha^{(1-\theta)}\beta^\theta>|\lambda_l|$. \vspace{0.1cm}

\textit{Claim 2: Let $\alpha=|\lambda_l|$, then $[\overline{l^2(E_l;\alpha)},\overline{l^2(E_l;\beta)}]_\theta=l^2(E_l;\alpha^{(1-\theta)}\beta^\theta)$.}

First, choose $\theta_1,\theta_2\in (0,1)$ such that $\theta_1+\theta_2-\theta_1\theta_2=\theta$, we can see
\[\begin{aligned}
{[}\overline{l^2(E_l;\alpha)},\overline{l^2(E_l;\beta)}{]}_\theta&=[l^2(E_l,A;\alpha),l^2(E_l;\beta)]_\theta\\
&=\big[[l^2(E_l,A;\alpha),l^2(E_l;\beta)]_{\theta_1},l^2(E_l;\beta)\big]_{\theta_2}\\
&\subset [l^2(E_l,A;\alpha^{(1-\theta_1)}\beta^{\theta_1}),l^2(E_l;\beta)]_{\theta_2}\\
&=[l^2(E_l;\alpha^{(1-\theta_1)}\beta^{\theta_1})\oplus \mathcal{S}_A(E_l),l^2(E_l;\beta)\oplus\{0\}]_{\theta_2}=l^2(E_l;\alpha^{(1-\theta)}\beta^\theta),
\end{aligned}\]
where we use Lemma \ref{lemma36} in the first equality, use the fact that $l^2(E_l;\beta)\subset l^2(E_l,A;\beta)$ in the third inequality, use Lemma \ref{lemma35} (b) in the fourth equality, and use Lemma \ref{lemma51} in the last equality. On the other hand, we also have 
\[ l^2(E_l;\alpha^{1-\theta}\beta^\theta)=[l^2(E_l;\alpha),l^2(E_l;\beta)]_\theta\subset [\overline{l^2(E_l;\alpha)},\overline{l^2(E_l;\beta)}]_{\theta}.\]
Combining the above two embedding relationships, noticing that both of them are continuous, we get the claim.\vspace{0.1cm}

Now we return to prove the lemma. We need to consider four different cases, based on whether $\alpha,\beta\in \{|\lambda_l|\}_{l=0}^\infty$. 

For the extreme case  $\alpha=|\lambda_l|,\beta=|\lambda_{l'}|$ for some $l'>l$, we devide the space $X$ into three pieces {$X=\hat{E}_l\oplus \hat{E}_{l'}\oplus Y$}, such that $\sigma(A,Y)=\sigma(A,X)\setminus {|\{z:|z|=|\lambda_l| \text{ or }|\lambda_{l'}|\}}$. By {using Lemma \ref{lemma36}}, we see that  
{\[\begin{cases}
\overline{l^2(X;\alpha)}=\overline{l^2(\hat{E}_l;\alpha)}\oplus \overline{l^2(\hat{E}_{l'};\alpha)}\oplus l^2(Y;\alpha),\\
\overline{l^2(X;\beta)}=\overline{l^2(\hat{E}_l;\beta)}\oplus \overline{l^2(\hat{E}_{l'};\beta)}\oplus l^2(Y;\beta).
\end{cases}
\]}
The assertion then follows by using Lemma \ref{lemma51} (a) and the two claims. The other three cases can be proved similarly. \hfill$\square$ \vspace{0.2cm}

The following are main theorems of this section.
\begin{theorem}\label{thm54}
	The Sobolev spaces $H^\sigma(K)$ are stable under complex interpolation.
\end{theorem}
\textit{Proof.} For $\sigma>\sigma'\geq0$, we  see that $$[H^\sigma(K),H^{\sigma'}(K)]_\theta=H^{(1-\theta)\sigma+\theta'\sigma}(K), \forall \theta\in[0,1]$$
is an easy consequence of Theorem \ref{thm47} (b), by using Lemma \ref{lemma51} (a), Lemma \ref{lemma52} and the fact that $l^2(X,A;\alpha)$ is stable under complex interpolation. \hfill$\square$\vspace{0.2cm}

The interpolation result for $H^\sigma_0(K)$ is somewhat complicated.  The result will coincide with the $\mathbb{R}^n$ case. Interested readers may read \cite{Lion} for an interpolation theorem for $H^\sigma_0(\Omega)$ with $\Omega\subset \mathbb{R}^n$.

{
\begin{definition}\label{def55}
	For $\sigma\geq 0$, define $$H^\sigma_{00}(K)=\big\{f\in H^{\sigma}(K):f\cdot\rho^{-\sigma}\in L^2(K)\big\},$$
	where $\rho(x)=R(x,V_0)^{1/2+d_H/2}$ on $K$ with $R(\cdot,\cdot)$ being the effective resistance metric. For each $f\in H^\sigma_{00}(K)$, we assign the norm
	\[\|f\|_{H^\sigma_{00}(K)}=\|f\|_{H^{\sigma}(K)}+\|f\rho^{-\sigma}\|_{L^2(K)}.\]
\end{definition}

{\noindent\textbf{Remark 1.} $H_{00}^\sigma(K)$ are natural analogs of the \textit{Lions-Magenes spaces}, see \cite{Lion}.}

\noindent\textbf{Remark 2.} One can easily see that for all $2k\geq\sigma$ with $k\in\mathbb{N}$, we have
\[H^\sigma_{00}(K)=\ker_\sigma\mathcal{T}_{\mathcal{H}_{k-1}}\oplus\big(\oplus_{\omega\in \mathcal{P}}\mathcal{R}_{\mathcal{H}_{k-1},\omega}\big(l^2(\mathcal{H}_{k-1};r_w^{\sigma/2}\mu_w^{(\sigma-1)/2})\big)\big).\]
In addition, for $f=f_0+\sum_{\omega\in \mathcal{P}}\mathcal{R}_{\mathcal{H}_{k-1},\omega}(\bm{h}_\omega) \in H_{00}^\sigma(K)$ with $f_0\in \ker_\sigma\mathcal{T}_{\mathcal{H}_{k-1}}$ and 
	$\bm{h}_\omega\in l^2(\mathcal{H}_{k-1};r_w^{\sigma/2}\mu_w^{(\sigma-1)/2})$, we have 
\[\|f\|_{H^\sigma_{00}(K)}\asymp\|f_0\|_{H^\sigma(K)}+\sum_{\omega\in\mathcal{P}}\|\bm{h}_\omega\|_{l^2(\mathcal{H}_{k-1};r_w^{\sigma/2}\mu_w^{(\sigma-1)/2})}.\]
}

We have the following interpolation theorems.

\begin{theorem}\label{thm56}
	The spaces $H_{00}^\sigma(K)$ are stable under complex interpolation.
\end{theorem}
\textit{Proof.} The theorem is a consequence of Remark 2, by a same method as Theorem \ref{thm54}.\hfill$\square$

\begin{theorem}\label{thm57}
	Let $\sigma>\sigma'\geq 0$, then $[H^\sigma_0(K),H^{\sigma'}_0(K)]_\theta=H^{(1-\theta)\sigma+\theta\sigma'}_{00}(K), \forall\theta\in(0,1).$ In particular, $[H^\sigma_0(K),H^{\sigma'}_0(K)]_\theta=H^{(1-\theta)\sigma+\theta\sigma'}_0(K)$ if and only if $r_w^{\sigma_\theta/2}\mu_w^{(\sigma_\theta-1)/2}\notin \{|\lambda_{l,w}|\}_{l=0}^\infty$  for any $\omega=\tau\dot{w}\in \mathcal{P}$, with $\sigma_\theta=(1-\theta)\sigma+\theta\sigma'$.
\end{theorem} 

\textit{Proof.} The first assertion follows from Theorem \ref{thm410}, Lemma \ref{lemma51} and \ref{lemma53}. The second assertion is a consequence of Lemma \ref{lemma36}. \hfill$\square$

\section{Interpolation of $H^\sigma(K)$: $\sigma\in\mathbb{R}$}
In this section, we will fulfill the definition of Sobolev spaces $H^\sigma(K)$ to negative orders, and study the associated interpolation theorem. Readers may read Lions and Magenes's monograph \cite{Lion} for classical theorems on bounded domains in $\mathbb{R}^n$. Part of the idea in this section is inspired by \cite{Lion}.

\begin{definition}\label{def61}
	For $\sigma\geq 0$, we define $H^{-\sigma}(K)=\big(H^\sigma_0(K)\big)'$, with the identification $H^0(K)=(H^0(K))'\subset H^{-\sigma}(K)$, noticing that $H^0(K)=H^0_0(K)=L^2(K)$. 
\end{definition}

\noindent\textbf{Remark 1.} In the above definition, we naturally embed the space $H^0(K)$ into $H^{-\sigma}(K)$. More concretely, for each $f\in H^0(K)$, we can correspond it with a linear functional $\varphi_f\in \big(H^\sigma_0(K)\big)'$ by the formula $\varphi_f(g)=\int_{K} \overline{f(x)}g(x)d\mu(x)=<g,f>_{L^2(K)}, \forall g\in H_0^\sigma(K)$. As a consequence, we always have 
\[H^\sigma(K)\subset H^{\sigma'}(K),\quad \forall \infty>\sigma>\sigma'>-\infty.\]

\noindent\textbf{Remark 2.} We also embed $H^0(K)$ into $\big(H^\sigma_{00}(K)\big)'$ in a same way. Notice that  $H^{-\sigma}(K)=\big(H^\sigma_{00}(K)\big)'$ if and only if $r_w^{\sigma/2}\mu_w^{(\sigma-1)/2}\notin\{|\lambda_{l,w}|\}_{l\geq 0}$ for any $\omega=\tau\dot{w}\in \mathcal{P}$.\vspace{0.1cm}

The following interpolation theorem is the main result in this section, which is a perfect analogue to the classical theorem.

\begin{theorem}\label{thm62}
	Let $-\infty<\sigma'<\sigma<+\infty$, $0<\theta<1$ and $\sigma_\theta=(1-\theta)\sigma+\theta\sigma'$. We have 
	\[[H^{\sigma}(K),H^{\sigma'}(K)]_{\theta}=\begin{cases}
	H^{\sigma_\theta}(K),\text{ if }\sigma_\theta\geq 0,\\
	\big(H_{00}^{-\sigma_\theta}(K)\big)',\text{ if }\sigma_\theta<0.
	\end{cases}\] 
	In particular, in case of $\sigma_\theta<0$, 
	$$[H^{\sigma}(K),H^{\sigma'}(K)]_{\theta}=H^{\sigma_\theta}(K)$$
    if and only if $r_w^{-\sigma_\theta/2}\mu_w^{-(\sigma_\theta+1)/2}\notin\{|\lambda_{l,w}|\}_{l\geq 0}$ for any $\omega=\tau\dot{w}\in \mathcal{P}$. 
\end{theorem}

In the rest of this section, we devote to prove Theorem \ref{thm62}.

\subsection{Lemmas}
We collect some lemmas first. For convenience, we let $Z$ be a Hilbert space with inner product $<\cdot,\cdot>_Z$. Let $Z_1\subset Z$ be a Banach space which is dense and continuously embedded in $Z$. 

Define $Z_{-1}$ as the dual space of $Z_1$, with the embedding $Z\subset Z_{-1}$ by 
\begin{equation}\label{eqn65}
z\to \varphi_z(\cdot)=<\cdot,z>_Z\in Z_{-1}.
\end{equation} 
By this, we have the relation $Z_1\subset Z\subset Z_{-1}$. In addition, we have the following lemma due to Lions and Magenes \cite{Lion} ({Proposition 2.1.}).

\begin{lemma}\label{lemma63}
	$[Z_{-1},Z_1]_{1/2}=Z$.
\end{lemma}

The next lemma will play a key role in the proof of Theorem \ref{thm62}.
\begin{lemma}\label{lemma64}
	Let $Z_1^{(0)}\subset Z_1$ be a closed subspace of $Z_1$ and suppose $Z_1^{(0)}$ is dense in $Z$. Define $Z_{-1}^{(0)}$ in a same way as $Z_{-1}$. If there is a map $L\in \mathcal{L}(Z,Z)\cap \mathcal{L}(Z_1,Z_1)$ such that 
	$$\begin{cases}
	L+id\in \mathcal{L}(Z,Z)\cap \mathcal{L}(Z_1,Z_1^{(0)}),\\
	L^*-id\in \mathcal{L}(Z,Z)\cap \mathcal{L}(Z_1,Z_1^{(0)}),
	\end{cases}$$
	where $id$ is the identity map and $L^*$ is the adjoint operator of $L$ with respect to $Z$, then 
	$$[Z^{(0)}_{-1},Z_1]_{1/2}=Z.$$
\end{lemma}

\textit{Proof.}  Let $\tilde{Z}={Z\times Z}$ with inner product $<\tilde{z},\tilde{z}'>_{\tilde{Z}}=<z_1,z'_1>_{Z}+<z_2,z_2'>_{Z}$, where $\tilde{z}=(z_1,z_2)$ and $\tilde{z}'=(z_1',z_2')$. Define $\tilde{Z}_1=\{(z_1,z_2)\in{Z}_1\times {Z}_1:z_1+z_2=Z_1^{(0)}\}$ with norm $\|(z_1,z_2)\|_{\tilde{Z}_1}=\|z_1\|_{Z_1}+\|z_2\|_{Z_1}$. Like we have done before, we define $\tilde{Z}_{-1}$ with $\tilde{Z}$ naturally embedded in {$\tilde{Z}_{-1}$}, noticing here that $\tilde{Z}_1$ is dense in $\tilde{Z}$ by the assumption on $Z_1^{(0)}$. 

Now we define the extension map $E\in \mathcal{L}(Z,\tilde{Z})\cap \mathcal{L}(Z_1,\tilde{Z}_{1})$ as follows
\[E(z)=(z,Lz), \quad\forall z\in Z.\]
The map $E$ can be naturally extended to be $E\in \mathcal{L}(Z^{(0)}_{-1},\tilde{Z}_{-1})$ with the  formula
\[E\varphi(\tilde{z})=\varphi(E^*{\tilde{z}}),\]
for any $\varphi\in Z^{(0)}_{-1}$ and $\tilde{z}\in\tilde{Z}_1$, noticing that $E^*\tilde{z}=E^*(z_1,z_2)=z_1+L^*z_2=z_1+z_2+(L^*-id)z_2\in Z^{(0)}_{1}$. Therefore, we get $E\in \mathcal{L}(Z^{(0)}_{-1},\tilde{Z}_{-1})\cap \mathcal{L}(Z_1,\tilde{Z}_1)$. As a consequence, $E\in \mathcal{L}([Z^{(0)}_{-1},Z_1]_{1/2},\tilde{Z})$ by using complex interpolation and Lemma \ref{lemma63}.

We also define a restriction map $R:\tilde{Z}_{-1}\to Z^{(0)}_{-1}$ by the formula
\[(R\tilde{\varphi})(z)=\tilde{\varphi}(z,0),\quad\forall \tilde{\varphi}\in \tilde{Z}_{-1} \text{ and }z\in Z^{(0)}_1.\]
It is then easy to see that $RE$ is the identity map from $Z^{(0)}_{-1}$ to $Z^{(0)}_{-1}$. In addition, we have $R(\tilde{Z})=Z$. Thus we get 
\[[Z^{(0)}_{-1},Z_1]_{1/2}=RE([Z^{(0)}_{-1},Z_1]_{1/2})\subset R(\tilde{Z})=Z.\]
On the other hand, we have $Z=[Z^{(0)}_{-1},Z^{(0)}_1]_{1/2}\subset[Z^{(0)}_{-1},Z_1]_{1/2}$ by using Lemma \ref{lemma63}. This finishes the proof. \hfill$\square$

\subsection{A decomposition by projection}
{Lemma \ref{lemma64} provides the strategy of the proof. Nevertheless, we need to overcome the difficulty that multiplication does not preserve smoothness in the fractal case \cite{BST}. In this part, for $k\in\mathbb{N}$, we focus on constructing a subspace $S_{k,w}$ in $L^2(K)$, such that the projections of functions in $H^{2k}(K)$ on ${S_{k,w}}$ maintain the smoothness, which will give us a new decomposition of the space $H^{2k}(K)$. This will play the role of multiplication by smooth bump functions.

	\begin{lemma}\label{lemma65} For $\omega=\dot{w}\in \mathcal{P}$, assuming $F_wK\cap V_0=\{\pi(\omega)\}$ without loss of generality, there is a linear map $\tilde{R}_{k,w}:\mathcal{H}_{k-1}\to dom\Delta^\infty$ such that $A_w\tilde{R}_{k,w}(h)=A_wh$, $P_{\tilde{R}_{k,w}(\mathcal{H}_{k-1})}(h)=\tilde{R}_{k,w}(h)$ and $\tilde{R}_{k,w}(h)$ is supported away from $V_0\setminus \{\pi(\omega)\}$.
   \end{lemma}	

\textit{Proof}. To achieve this, we choose a basis $\{h_1,h_2,\cdots,h_n\}$ of $\mathcal{H}_{k-1}$, and denote 
	\[a_{ij}=<h_i,h_j>_{L^2(F_wK)},\quad 1\leq i,j\leq n.\]
It is clear that we can find $\tilde{h}_i\in dom\Delta^\infty,1\leq i\leq n$, such that $A_w\tilde{h}_i=A_wh_i$, the support of $\tilde{h}_i$ is a small neighbourhood of $F_wK$, and $<\tilde{h}_i,h_j>_{L^2(K)}=a_{ij}$. In addition, we can assume that 
 \[<\tilde{h}_i,\tilde{h}_j>_{L^2(K)}=\varepsilon_{ij}+a_{ij},\quad  1\leq i,j\leq n\]
with $\varepsilon=\max_{i,j}\{|\varepsilon_{ij}|\}$  small enough so that we can find $f_i\in dom\Delta^\infty$ supported in some compact subsets of $K\setminus F_wK$ away from the boundary, satisfying 
\[\begin{cases}
<f_i,\tilde{h}_j>_{L^2(K)}=0, &\forall 1\leq i,j\leq n,\\
<f_i,f_j>_{L^2(K)}=\delta_{ij}\varepsilon, &\forall 1\leq i,j\leq n,\\
<f_i,h_j>_{L^2(K)}=\varepsilon_{ij}+\delta_{ij}\varepsilon, &\forall 1\leq i,j\leq n.
\end{cases}\]
Set $\tilde{R}_{k,w}(h_i)=\tilde{h}_i+f_i$, and extend $\tilde{R}_{k,w}$ to be the linear map $\mathcal{H}_{k-1}\to dom\Delta^\infty$. One can then check that \[<h_i,\tilde{R}_{k,w}(h_j)>_{L^2(K)}=<\tilde{R}_{k,w}(h_i),\tilde{R}_{k,w}(h_j)>_{L^2(K)}=a_{ij}+\varepsilon_{ij}+\delta_{ij}\varepsilon,\]
and thus $P_{\tilde{R}_{k,w}(\mathcal{H}_{k-1})}(h_i)=\tilde{R}_{k,w}(h_i)$ for any $1\leq i\leq n$. The lemma follows immediately. \hfill$\square$

\begin{definition}\label{def63} Let $\omega=\dot{w}\in\mathcal{P}$, $k\geq 1$ be an integer.

(a). Write $f_{w,h}$ for $\tilde{R}_{k,w}(h)$ for short. We omit $k$ since $\tilde{R}_{k,w}$ can defined consistently for different $k$'s.
	
(b). Let $S_{k,w}$ be the subspace of $L^2(K)$ spanned by the functions $\{f_{w,h}\circ F_w^{-n}:h\in \mathcal{H}_{k-1},n\geq 0\}$.
\end{definition}

We have the following theorem. 
\begin{theorem}\label{prop66} Let $k\geq 1$ be an integer, $f\in H^{2k}(K)$, and $\omega=\tau\dot{w}\in \mathcal{P}$.
	
	(a). If $A_\tau f\bot S_{k,w}$ in $L^2(K)$, we have $T_{\omega}^{(2k)}f=0$.
	
	(b). $P_{S_{k,w}}A_\tau f\in H^{2k}(K)$, and $f-\mu_\tau A_\tau^*P_{S_{k,w}}A_\tau f\in H^{2k}(K)$, where $A^*_\tau$ is the adjoint operator of $A_\tau$ {in $L^2(K)$, which can be expressed by $A_\tau^*g=\mu^{-1}_{\tau}g\circ F^{-1}_\tau, \forall g\in L^2(K)$.}
\end{theorem}
\textit{Proof.} (a). Let $h=T_{\omega}^{(2k)}f$ and $\tilde{f}=A_\tau f-h$. Then we have 
\[\begin{aligned}
<f_{w,A_w^nh}\circ F_w^{-n}\circ F_\tau^{-1},f>_{L^2(K)}&=\mu_w^n\mu_\tau<f_{w,A_w^nh},A_w^nA_\tau f>_{L^2(K)}\\
&=\mu_w^n\mu_\tau\big(<f_{w,A_w^nh},A_w^nh>_{L^2(K)}+<f_{w,A_w^nh},A_w^n\tilde{f}>_{L^2(K)}\big)\\
&=\mu_w^n\mu_\tau \big(\|f_{w,A_w^nh}\|_{L^2(K)}^2+o(\lambda_{l_w(2k),\omega}^n)\|f_{w,A_w^nh}\|_{L^2(K)}\big).
\end{aligned}\]
Thus the left side equals $0$ for any $n\geq 0$ only if $T_{\omega}^{(2k)}f=h=0$.

(b). Without loss of generality, we consider the case that $\omega=\dot{w}$. For each $f\in H^{2k}(K)$, we will construct a series $\sum_{n=0}^\infty f^{(n)}$ converging in $H^{2k}(K)$, where each $f^{(n)}$ takes the form $f^{(n)}=f_{w,h}\circ F_w^{-n}$ for some $h\in\mathcal{H}_{k-1}$, so that $P_{S_{k,w}}f=\sum_{n=0}^\infty f^{(n)}$.

First, we look at some special functions. Let $f\in L^2(K)$ such that $A_w^lf\in\mathcal{H}_{k-1}$ for some $l\geq 0$. Denote $S^{(l)}_{k,w}=\{f_{w,h}\circ F_w^{-n}:h\in \mathcal{H}_{k-1},0\leq n\leq l\}$, and write $P_{S_{k,w}^{(l)}}f=\sum_{n=0}^{l} f^{(n)}$. Clearly, $g=f-\sum_{n=0}^{l-1} f^{(n)}$ is $k$-multiharmonic in $F_w^{l}K$, and so $f^{(l)}=f_{w,A_w^{l}g}\circ F_w^{-l}$ by Lemma \ref{lemma65}. As a consequence, we have $f-P_{S_{k,w}^{(l)}}f=0$ on $F_{w}^{l+1}K$, which shows that $f-P_{S_{k,w}^{(l)}}f\bot S_{k,w}$. By this observation, we have the following construction.}\vspace{0.1cm}

\textit{Step 1: For any $f\in L^2(K)$ such that $A_w^lf\in\mathcal{H}_{k-1}$ for some $l\geq 0$, we can write $P_{S_{k,w}}f=\sum_{n=0}^{l} f^{(n)}$, where each $f^{(n)}$ takes the form $f^{(n)}=f_{w,h}\circ F_w^{-n}$ for some $h\in\mathcal{H}_{k-1}$.}\vspace{0.1cm}

\textit{Step 2: For any $f\in L^2(K)$ such that $A_w^lf\in\mathcal{H}_{k-1}$ for some $l\geq 0$, we have by induction
\[\|f^{(n)}\|_{L^\infty(K)}\lesssim\|\sum_{m=0}^{n}f^{(m)}\|_{L^2(F^{n}_wK\setminus F_{w}^{n+1}K)}+\sum_{m=0}^{n-1}\|f^{(m)}\|_{L^\infty(K)}\lesssim {2^n} \|f\|_{L^2(K)},\text{ for } n\geq 0.
\] 
So we can continuously extend the definition of $f^{(n)},n\geq 0$ to general functions $f$ in $L^2(K)$. }\vspace{0.1cm} 

We have some observations {on} the sequence $\{f^{(n)}\}_{n\geq 0}$. \vspace{0.1cm}

\textit{Observation 1: For any $f\in L^2(K)$ and $n\geq 1$, $A_w^{n-1}f^{(n)}=(A_w^{n-1}f)^{(1)}$.} 

\textit{Proof of Observation 1.}  Only need to consider the case that $A_w^lf\in \mathcal{H}_{k-1}$ for some $l\geq 0$. {Let $g=f-\sum_{m=0}^{n-2}f^{(m)}$. Then we have 
\[A_w^{n-1}P_{S_{k,w}}g=A_w^{n-1}P_{S_{k,w}^{(n-1)+}}g=P_{S_{k,w}}\big(A_w^{n-1}g\big),\]
where $S^{(n-1)+}_{k,w}=\{f_{w,h}\circ F_w^{-m}:h\in \mathcal{H}_{k-1},m\geq n-1\}$.
So we have $A_w^{n-1}f^{(n)}=A_w^{n-1}g^{(n)}=(A_w^{n-1}g)^{(1)}$. On the other hand, we have $(A_w^{n-1}f)^{(1)}=(A_w^{n-1}g)^{(1)}$ as $A_w^{n-1}(f-g)\in\mathcal{H}_{k-1}$.}\vspace{0.1cm} 

\textit{Observation 2: There a kernel $\psi\in L^\infty(K\times K)$ such that
$$A_w^{n-1}f^{(n)}(x)=\int_K \psi(x,y)\Delta^k \big(A_w^{n-1}f(y)\big)d\mu(y),$$
for any $f\in H^{2k}(K)$ and $n\geq 1$. }

\textit{Proof of Observation 2.} We only need to choose
$$
\psi(x,y)=(-1)^k\int_{K^{k-1}}G^{(1)}_{y_1}(x)G(y_1,y_2)\cdots G(y_{k-1},y)d\mu(y_1)\cdots d\mu(y_{k-1}),
$$
where $G_y(x)=G(x,y)$ is the Green's function. By Step 2, we immediately have $\psi\in L^\infty(K\times K)$. Since $h^{(1)}=0,\forall h\in \mathcal{H}_{k-1}$, we can easily see that $f^{(1)}(x)=\int_K \psi(x,y)\Delta^kf(y)d\mu(y)$. For $n\geq 2$, we use Observation 1.\vspace{0.1cm}

Now, by using Observation 2 and Lemma \ref{lemma38}, we can see that $\{f^{(n)}\}_{n\geq 0}\in l^2(L^\infty(K);r_w^k\mu_w^{k-1/2})$ for any $f\in H^{2k}(K)$. Then, a same proof as in Theorem \ref{thm45} shows that $\sum_{n=0}^{\infty} f^{(n)}$ converges in  $H^{2k}(K)$ with $\|\sum_{n=0}^{\infty} f^{(n)}\|_{H^{2k}(K)}\lesssim \|f\|_{H^{2k}(K)}$. On the other hand, we have $P_{S_{k,w}}f=\sum_{n=0}^\infty f^{(n)},\forall f\in H^{2k}(K)$ as desired. In fact, this is ture if $A_w^lf\in \mathcal{H}_{k-1}$ for some $l\geq 0$, and this kind of functions are dense in $H^{2k}(K)$. \hfill$\square$\vspace{0.1cm} 

\subsection{Proof of Theorem \ref{thm62}} {We return to prove  Theorem \ref{thm62}. As shown in Lemma \ref{lemma64}, the key is to construct the map `$L$'. Theorem \ref{prop66} will play a crucial role.}

\begin{lemma}\label{lemma67}
	For each $w\in W_*$ and $l\geq 0$, there is a polynomial $p_w$ such that $p_w(A_w)+id=0$ and $p_w(\mu_w^{-1}A_w^{-1})-id=0$ on $\tilde{E}_{l,w}$.
\end{lemma}
\textit{Proof.} Since $\tilde{E}_{l,w}$ is of finite dimensional, there are polynomials $p_1$ and $p_2$ such that $p_1(A_w)=0$ and $p_2(\mu_w^{-1}A_w^{-1})=0$ on $\tilde{E}_{l,w}$. The zeros of $p_1$ and $p_2$ can be disjoint, since they can  be just the eigenvalues of $A_w$ and $\mu_w^{-1}A_w^{-1}$ respectively. Then $p_1$ and $p_2$ are coprime polynomials, and thus there exist polynomials $r_1$ and $r_2$ such that $r_1p_1-r_2p_2=1$. Then the polynomial  $p_w=r_1p_1+r_2p_2$ will satisfy the requirement of the lemma.\hfill$\square$

\begin{definition}\label{def68}
(a). Let $k\geq 1$, $\omega=\tau\dot{w}\in\mathcal{P}$ and $l=l_{\omega(2k)}$.  Take $p_w$ as in Lemma \ref{lemma67}.  Define $L^{(k)}_\omega=\mu_{\tau}A_\tau^* P_{S_{k,w}}p_w(A_w)A_\tau,$ where $A_\tau^*$ is the adjoint operator of $A_\tau$. 

(b). Define $L^{(k)}=\sum_{\omega\in \mathcal{P}} L_\omega^{(k)}$.
\end{definition}

We have the following Proposition.
\begin{proposition}\label{prop69}
	$L^{(k)}+id\in \mathcal{L}\big(H^{0}(K),H^{0}(K)\big)\cap \mathcal{L}\big(H^{2k}(K),H^{2k}_0(K)\big)$ and $(L^{(k)})^*-id\in \mathcal{L}\big(H^{0}(K),H^{0}(K)\big)\cap \mathcal{L}\big(H^{2k}(K),H^{2k}_0(K)\big)$, where $id$ is the identity map.
\end{proposition}
\textit{Proof}. Let $\omega=\tau\dot{w}\in\mathcal{P}$. For any $0\leq s<\infty$ and $f\in H^{2k}(K)$, by using Theorem \ref{prop66}, we always have  $\mu_\tau A_\tau^* P_{S_{k,w}}A_w^s A_\tau f\in H^{2k}(K)$, has $0$ tangent on $\mathcal{P}\setminus \{\pi(\omega)\}$, and 
\[T_\omega^{(2k)}(\mu_\tau A_\tau^* P_{S_{k,w}}A_w^sA_\tau f)=T_\omega^{(2k)}(\mu_\tau A_\tau^* A_w^sA_\tau f)=A_w^s\big(T_\omega^{(2k)}(f)\big).\]
As a consequence, we see that $L^{(k)}_\omega f$ has $0$ tangent at $\mathcal{P}\setminus \{\pi(\omega)\}$, and $T_\omega^{(2k)}(L^{(k)}_\omega f)=p_w(A_w)T_\omega^{(2k)}(f)$. By Definition \ref{def68} and Lemma \ref{lemma67}, we conclude that $(L^{(k)}+id)f\in H^{2k}_0(K)$ using the characterization of $H^{2k}_0(K)$ in Theorem \ref{thm42}.

To show the other half, we need some observations.

\textbf{1)} $\big(L^{(k)}_\omega\big)^*=\mu_\tau A_\tau^* p_w(A_w^*)P_{S_{k,w}}A_\tau,$ and $\big(L^{(k)}\big)^*=\sum_{\omega\in\mathcal{P}} \big(L^{(k)}_\omega\big)^*$.

\textbf{2)} For $f\in H^{2k}(K)$, $T_\omega^{(2k)}(A_w^*f)=\mu_w^{-1}A_w^{-1}\big(T_\omega^{(2k)}(f)\big)$, noticing that $A_w^*(f)=\mu_w^{-1} f\circ F_w^{-1}$.  

%%\textbf{3)} Since $m$ is large enough for each $\omega$, we have $\big(L^{(k)}_\omega \big)^*f\in H^{2k}(K)$ for any $f\in H^{2k}(K)$.

The rest of the proof is similar to that of the first part. \hfill$\square$\vspace{0.1cm}

\textit{Proof of Theorem \ref{thm62}.} By using Lemma \ref{lemma64} and Proposition \ref{prop69}, we see that 
\[[H^{-2k}(K),H^{2k}(K)]_{1/2}=H^0(K).\] 
Thus for $-\infty<\sigma_1<0<\sigma_2<+\infty$ and $\theta=\frac{\sigma_1}{\sigma_1-\sigma_2}$, by using Theorem \ref{thm54} and \ref{thm57}, we get
\[\big[\big(H_{00}^{-\sigma_1}(K)\big)',H^{\sigma_2}(K)\big]_\theta=H^0(K).\]
As a consequence, we then have $$[H^{\sigma_1}(K),H^{\sigma_2}(K)]_\theta=H^0(K),$$
because 
$$H^0(K)=[H^{\sigma_1}(K),H^{\sigma_2}_0(K)]_\theta\subset [H^{\sigma_1}(K),H^{\sigma_2}(K)]_\theta\subset \big[\big(H_{00}^{-\sigma_1}(K)\big)',H^{\sigma_2}(K)\big]_\theta=H^0(K).$$
Then the theorem follows immediately.\hfill$\square$

\section{Examples}
In this section, we present some concrete examples.

\subsection{D3-symmetric fractals}
Tangents on D3-symmetric p.c.f. self-similar sets have been studied in detail in \cite{s3} on the domain of $\Delta^k$, with some related studies in \cite{cq1, cq4, distribution}. 

\begin{figure}[htp]
	\includegraphics[height=3.6cm]{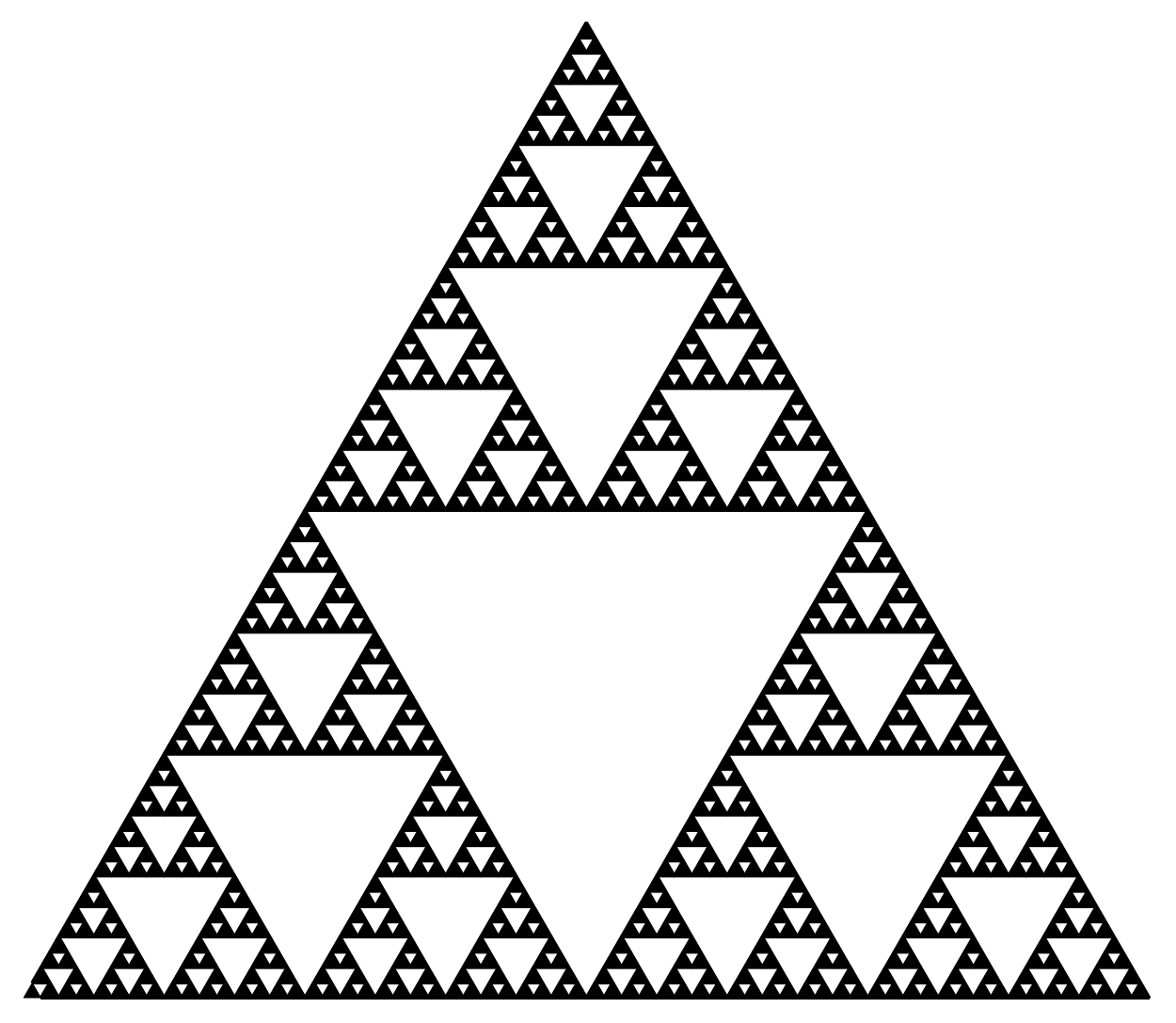}\quad
	\includegraphics[height=3.6cm]{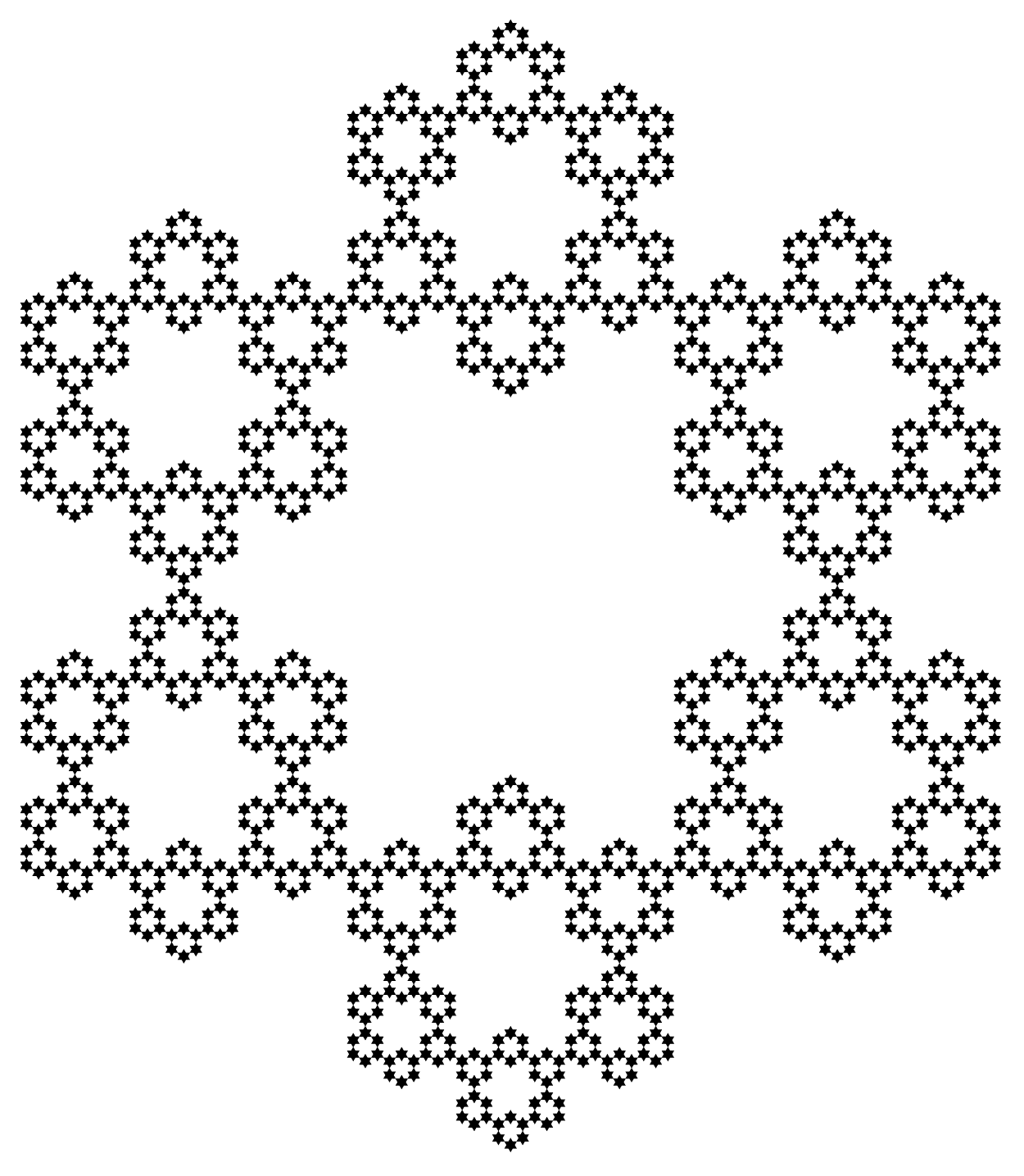}\quad
	\includegraphics[height=3.6cm]{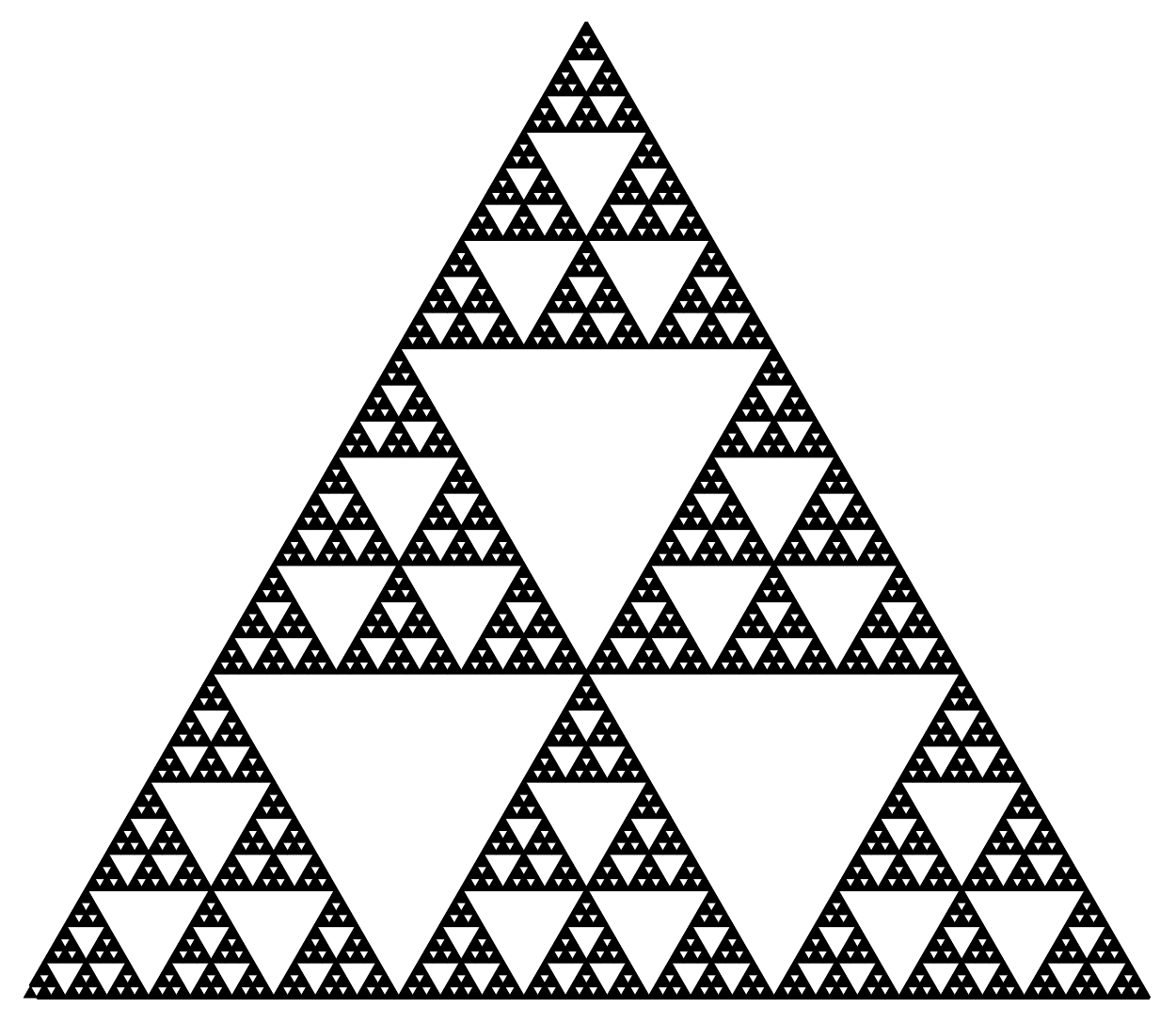}
	\caption{Examples of D3-symmetric p.c.f. self-similar sets: the Sierpinski gasket,  the Hexagasket and the level-$3$ Sierpinski gasket.}\label{fig1}
\end{figure}

More precisely, let's look at a p.c.f. self-similar set $K$ with exactly three boundary points $V_0=\{p_1,p_2,p_3\}$ {such that} $\pi^{-1}(p_i)=\dot{i}$. Assume that there exists a group $\mathcal{G}$ of homeomorphisms of $K$ isomorphic to the D3-symmetric group that acts as permutations on $V_0$, and $\mathcal{G}$ preserves the harmonic structure and the self-similar measure of $K$. See Figure \ref{fig1} for examples.

For fixed $i\in\{1,2,3\}$, let $h_T$ be the antisymmetric harmonic function with the boundary values $h_T(p_i)=0,h_T(p_{i+1})=1,h_T(p_{i+2})=-1$, where we use the cyclic notation $p_4=p_1$. Then it is easy to see that,
\[\{\lambda_{l,i}\}_{l\geq 0}=\{r_i^n\mu_i^n,r_i^{n+1}\mu_i^n,\iota_ir_i^n\mu_i^n\}_{n\geq 0},\]
where $\iota_i$ is defined by the identity  $A_ih_T=\iota_ih_T$.

In convention, we define \textit{normal derivatives} and \textit{tangential derivatives} of functions at $p_i$ by the following pointwise formulas, if the limits exist,
\[\begin{cases}
\partial_n f(p_i)=\lim\limits_{n\to\infty} r_i^{-n}\big(2f(p_i)-f(F_i^np_{i+1})-f(F_i^np_{i+2})\big),\\
\partial_T f(p_i)=\lim\limits_{n\to\infty} \iota_i^{-n}\big(f(F_i^np_{i+1})-f(F_i^np_{i+2})\big).
\end{cases}\]
 Assuming (\textbf{A2}), by using Theorem \ref{thm314}, we can easily see the following result. 

\begin{theorem}
(a). For $\sigma>2n+2-\frac{d_S}{2}$, $n\in \mathbb{Z}^+$, $\partial_n\Delta^nf(p_i)$ is well defined,  $\forall f\in H^\sigma(K)$.

(b). For $\sigma>2n+\frac{2\log \iota_i}{(1+d_H)\log r_i}+\frac{d_S}{2}$, $n\in \mathbb{Z}^+$, $\partial_T\Delta^nf(p_i)$ is well defined, $\forall f\in H^\sigma(K)$.
\end{theorem} 

The following is an equivalent narration of Theorem \ref{thm42}.
\begin{theorem}
For $\sigma\geq 0$ and $f\in H^\sigma(K)$, we have $f\in H^\sigma_0(K)$ if and only if
$$\begin{cases}
\Delta^n f(p_i)=0,&\forall 0\leq n<\frac{\sigma}{2}-\frac{d_S}{4}\text{ and }i=1,2,3,\\
\partial_n\Delta^n f(p_i)=0,&\forall 0\leq n<\frac{\sigma}{2}+\frac{d_S}{4}-1\text{ and }i=1,2,3,\\
\partial_T\Delta^n f(p_i)=0,&\forall 0\leq n<\frac{\sigma}{2}-\frac{d_S}{4}-\frac{\log \iota_i}{(1+d_H)\log r_i}\text{ and }i=1,2,3.
\end{cases}$$
\end{theorem}

\subsection{The Vicsek set}
\begin{figure}[htp]
	\includegraphics[height=3.6cm]{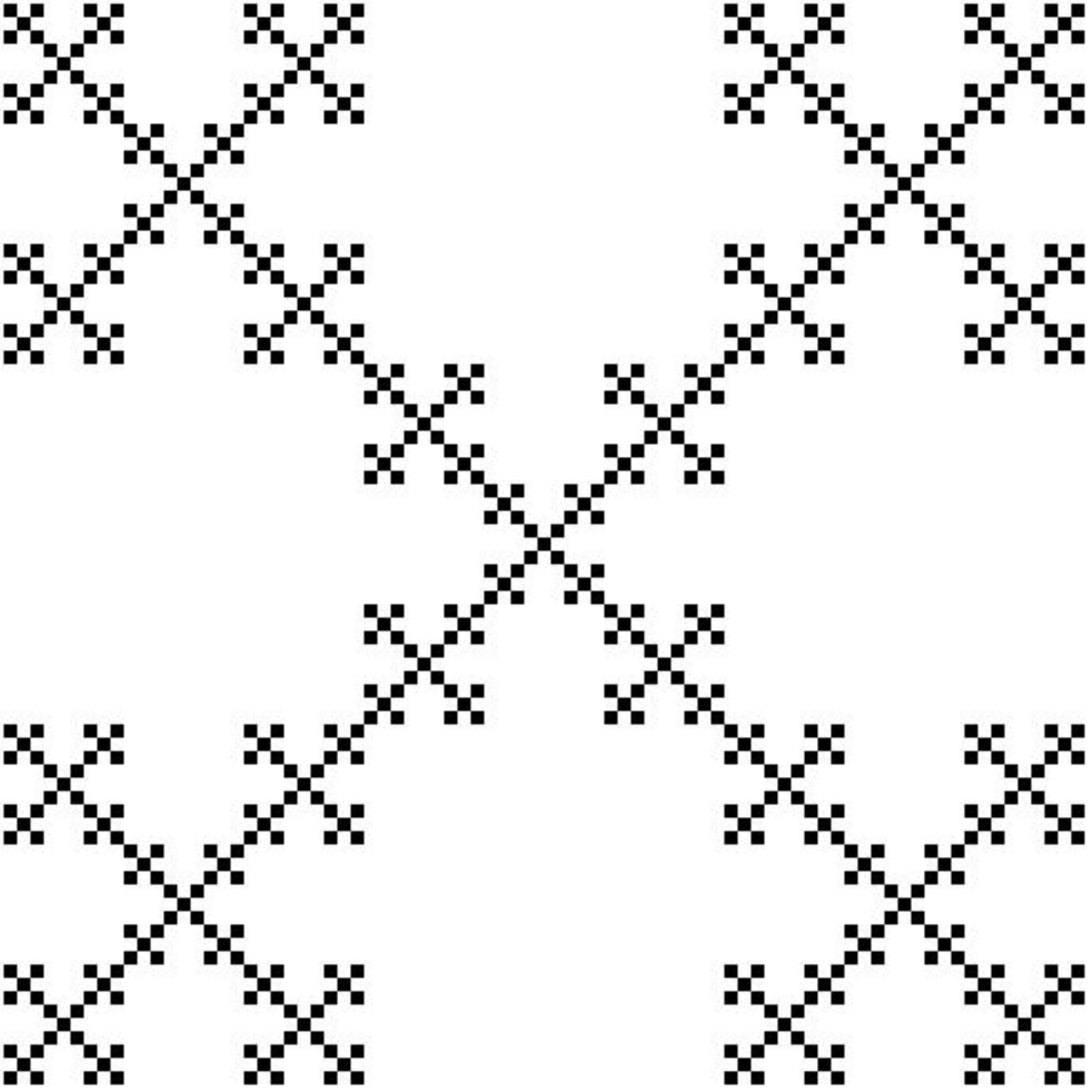}\quad
	\caption{The Vicsek set $\mathcal{V}$.}\label{fig2}
\end{figure}

Let $\{p_i\}_{i=1}^4$ be the four vertices of a unit square, and $p_5$ be the center of the square. The Vicsek set $\mathcal{V}$ (see Figure \ref{fig2}) is the attractor of the i.f.s. $\{F_i\}_{i=1}^5$, where
\[F_ix=\frac{1}{3}x+\frac{2}{3}p_i,\text{ for }i=1,2,3,4,5.\]

The boundary set of $\mathcal{V}$ is $V_0=\{p_i\}_{i=1}^4$ with $\pi^{-1}(p_i)=\dot{i}$. There is a unique S4-symmetric harmonic structure on $\mathcal{V}$, with 
\[r_i=\frac{1}{3},i=1,2,3,4,5,\]
and
\[\mathcal{E}_0(f,g)=\sum_{i\neq j} \big(f(p_i)-f(p_j)\big)\big(g(p_i)-g(p_j)\big),\quad \forall f,g\in l(V_0).\]
In addition, we take $\mu$  to be the canonical normalized Hausdoff measure on $\mathcal{V}$.

The Vicsek set $\mathcal{V}$ is an interesting example in that $\{\lambda_{l,i}\}_{l\geq 0}=\{15^{-n},3^{-1}\cdot 15^{-n}\}_{n\geq 0}$, with each $\lambda_{l,i}$ has a one dimensional generalized eigenspace. We have the following narration of Theorem \ref{thm42}.

\begin{theorem}
	For $\sigma\geq 0$ and $f\in H^\sigma(\mathcal{V})$, we have $f\in H^\sigma_0(\mathcal{V})$ if and only if
	$$\begin{cases}
	\Delta^n f(p_i)=0,&\forall 0\leq n<\frac{\sigma}{2}-\frac{d_S}{4}\text{ and }i=1,2,3,4,\\
	\partial_n\Delta^n f(p_i)=0,&\forall 0\leq n<\frac{\sigma}{2}+\frac{d_S}{4}-1\text{ and }i=1,2,3,4.
	\end{cases}$$
\end{theorem}

Furthermore, for $\sigma\geq 0$, write $H^\sigma_D(\mathcal{V})=(id-\Delta_D)^{-\sigma/2}L^2(\mathcal{V})$ and $H^\sigma_N(\mathcal{V})=(id-\Delta_N)^{-\sigma/2}L^2(\mathcal{V})$, where $\Delta_D$ and $\Delta_N$ are the Dirichlet and Neumann Laplacians. See \cite{cq3,s1} for a detailed discussion on these spaces. Then we have 

\begin{theorem}
	For $\sigma\geq 0$, $H_{00}^\sigma(\mathcal{V})=H^\sigma_D(\mathcal{V})\cap H^\sigma_N(\mathcal{V})$ with $\|f\|_{H_{00}^\sigma(\mathcal{V})}\asymp \|f\|_{H^\sigma_D(\mathcal{V})}+\|f\|_{H^\sigma_N(\mathcal{V})}$.
\end{theorem}
\textit{Proof.} Fix $k\geq \lceil \sigma/2\rceil$, we break the spaces $\mathcal{H}_{k-1}$ into two parts, $\mathcal{H}_{k-1}=X^{(i)}\oplus Y^{(i)}$ such that $\sigma(A_i;X^{(i)})=\{1,15^{-1},\cdots,15^{-k+1}\}$ and $\sigma(A_i;Y^{(i)})=\{3^{-1},45^{-1},\cdots,3^{-1}\cdot15^{-k+1}\}$. 

Using the notations in Definition \ref{def46}, one can check that 
\[
\begin{cases}
\begin{aligned}H^\sigma_D(\mathcal{V})=\ker_\sigma\mathcal{T}_{\mathcal{H}_{k-1}}&\oplus\Big(\oplus_{i=1}^4\mathcal{R}_{\mathcal{H}_{k-1},\dot{i}}\big(l^2(X^{(i)};r_i^{\sigma/2}\mu_i^{(\sigma-1)/2})\big)\Big)\\&\oplus\Big(\oplus_{i=1}^4\mathcal{R}_{\mathcal{H}_{k-1},\dot{i}}\big(l^2(Y^{(i)},A_i;r_i^{\sigma/2}\mu_i^{(\sigma-1)/2})\big)\Big),\end{aligned}\\
\begin{aligned}H^\sigma_N(\mathcal{V})=\ker_\sigma\mathcal{T}_{\mathcal{H}_{k-1}}&\oplus\Big(\oplus_{i=1}^4\mathcal{R}_{\mathcal{H}_{k-1},\dot{i}}\big(l^2(X^{(i)},A_i;r_i^{\sigma/2}\mu_i^{(\sigma-1)/2})\big)\Big)\\&\oplus\Big(\oplus_{i=1}^4\mathcal{R}_{\mathcal{H}_{k-1},\dot{i}}\big(l^2(Y^{(i)};r_i^{\sigma/2}\mu_i^{(\sigma-1)/2})\big)\Big).
\end{aligned}
\end{cases}
\]
The theorem follows immediately. \hfill$\square$
%\noindent{\bf Acknowledgment.} ***.

\section{On the assumption (\textbf{A1})}
Although a large amount of p.c.f. fractals, including all the examples in Section 7, satisfy (\textbf{A1}), there do exist counter examples. \vspace{0.1cm}

\begin{figure}[htp]
	\includegraphics[height=3.6cm]{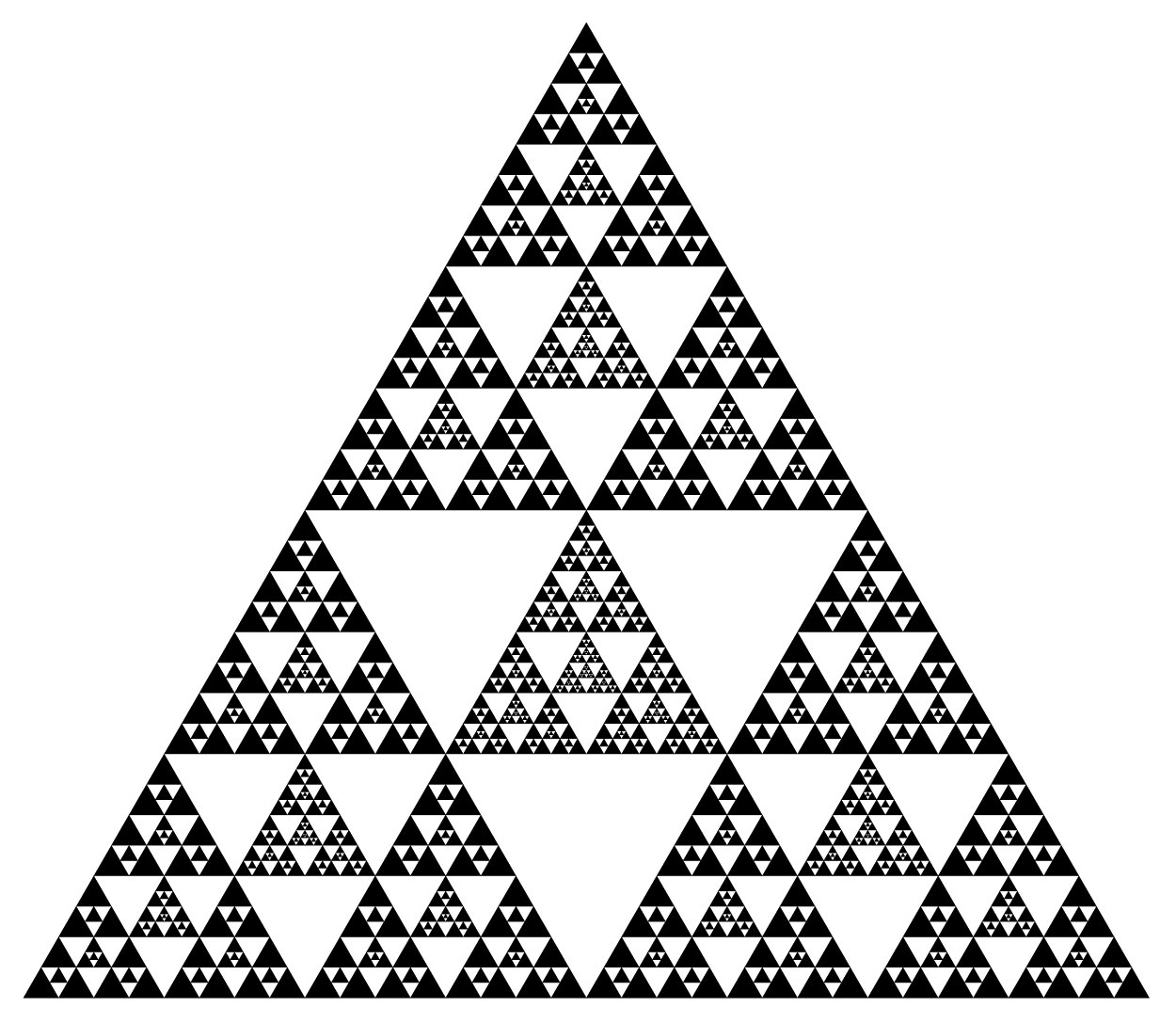}\quad
	\caption{The filled Sierpinski gasket $\mathcal{SG}_f$.}\label{fig3}
\end{figure}

\noindent\textbf{Example.} Let $\{p_1,p_2,p_3\}$ be the three vertices of a triangle, and $p_4=\frac{1}{3}\sum_{i=1}^3 p_i$  be the center. We define an i.f.s. $\{F_i\}_{i=1}^4$ by 
\[\begin{aligned}
F_i(x)&=\frac{1}{2}x+\frac{1}{2}p_i,  i=1,2,3,\\
F_4(x)&=\frac{1}{4}x+\frac{3}{4}p_4.
\end{aligned}\]
Call the  unique compact set, denoted by $\mathcal{SG}_f$, satisfying
$\SG_f=\bigcup_{i=1}^4 F_i(\SG_f)$, the \textit{filled Sierpinski gasket}. 
See Figure \ref{fig3}. 

One can check that $\mathcal{C}=\{1\dot{2},1\dot{3},2\dot{1},2\dot{3},3\dot{1},3\dot{2},4\dot{1},4\dot{2},4\dot{3},12\dot{3},13\dot{2},21\dot{3},23\dot{1},31\dot{2},32\dot{1}\}$,  $\mathcal{P}=\{\dot{1},\dot{2},\dot{3},1\dot{2},1\dot{3},2\dot{1},2\dot{3},3\dot{1},3\dot{2}\}$ and $V_0=\{p_1,p_2,p_3,F_1p_2,F_2p_3,F_3p_1\}$. One can see that $$\pi^{-1}(F_1p_2)=\{1\dot{2},2\dot{1}\}.$$
As a consequence, (\textbf{A1}) fails for $\SG_f$.\hfill$\square$\vspace{0.1cm}

Fortunately, \textbf{all the main theorems} in this paper, including Theorem \ref{thm314}, \ref{thm42}, \ref{thm54}, \ref{thm56}, \ref{thm57} and \ref{thm62}, are valid even if (\textbf{A1}) is not satisfied, as long as we assume (\textbf{A2}). Clearly, we did not use (\textbf{A1}) in Section 3, but Section 4 and Section 5 are somewhat delicate, where we need a precise description of the pretangents at the boundary. Below we briefly show the necessary materials in proving Theorem \ref{thm42}, \ref{thm54}, \ref{thm56}, \ref{thm57} and \ref{thm62} without using (\textbf{A1}). 

We need some new notations. Let $\omega=\dot{w}\in\mathcal{P}$ and  $k\in\mathbb{N}$.

\noindent\textbf{Notation 1:} \textit{Denote $\tilde{A}_w: \mathbb{C}^k\to\mathbb{C}^k$ by $(x_0,x_1,\cdots,x_{k-1})\to (x_0,r_w\mu_wx_1,\cdots, (r_w\mu_w)^{k-1}x_{k-1})$.}

\noindent\textbf{Notation 2:} \textit{{Denote ${\mathcal{H}}^{(w)}_{k-1}=\{h\in\mathcal{H}_{k-1}:\Delta^l h(\pi(\dot{w}))=0,\forall 0\leq l\leq k-1\}$.}}

With some effort, one can check the following claims. 

\noindent\textbf{Claim 1:} \textit{There is a natural isomorphism $ \mathbb{C}^k\times \mathcal{H}^{(w)}_{k-1}\to \mathcal{H}_{k-1}$, which gives us natural isomorphisms $I_{w}:l^2(\mathbb{C}^k,\tilde{A}_w;\alpha)\times l^2({\mathcal{H}}^{(w)}_{k-1},A_w;\alpha)\to l^2(\mathcal{H}_{k-1},A_w;\alpha)$.}

\noindent\textbf{Claim 2:} \textit{Let $\sigma\geq0$ and $k\geq \lceil \sigma/2\rceil$, $p\in V_0$, $\pi^{-1}(p)=\{\omega_1,\cdots,\omega_m\}=\{\tau_1\dot{w}_1,\cdots, \tau_m\dot{w}_m\}$ and 
	$$\begin{aligned}
	&X_{p,\sigma}=\big\{  \big(\{\bm{x}_n^{(1)}\}_{n\geq 0},\{\bm{x}_n^{(2)}\}_{n\geq 0},\cdots,\{\bm{x}_n^{(m)}\}_{n\geq 0}\big)\in {\prod_{i=1}^m} l^2(\mathbb{C}^k,\tilde{A}_{w_i};r_{w_i}^{\sigma/2}\mu_{w_i}^{(\sigma-1)/2}):\\&\lim\limits_{n\to\infty} (r_{\tau_i}\mu_{\tau_i})^{-l}(r_{w_i}\mu_{w_i})^{-nl}(\bm{x}^{(i)}_{n})_l=\lim\limits_{n\to\infty} (r_{\tau_j}\mu_{\tau_j})^{-l}(r_{w_j}\mu_{w_j})^{-nl}(\bm{x}^{(j)}_{n})_l,\forall i\neq j, 0\leq l\leq k-1\big\}.
	\end{aligned}$$ 
There is an isomorphism 
	$$I_p:X_{p,\sigma}\to l^2(\mathbb{C}^k,\tilde{A}_{w_1};r_{w_1}^{\sigma/2}\mu_{w_1}^{(\sigma-1)/2})\times \big(\prod_{i=2}^m l^2(\mathbb{C}^k;r_{w_i}^{\sigma/2}\mu_{w_i}^{(\sigma-1)/2})\big),$$
defined consistently for all $\sigma\neq 2l+\frac{d_S}{2}, 0\leq l\leq k-1$.}

The key steps of constructing $I_p$ is: first we pick a {nondecreasing sequence} $\{\ell_n\}_{n\geq0}$ such that $r_{w_1}^{\ell_n}\asymp r_{w_j}^n$, and hence by (\textbf{A2}) $\mu_{w_1}^{\ell_n}\asymp \mu_{w_j}^n$; next we define a sequence $\{\hat{\bm{x}}_n^{(j)}\}_{n\geq 0}$ by $$A_{\tau_j}^{-1}\big(\tilde{A}_{w_j}^{-n-1}\hat{\bm{x}}^{(j)}_{n+1}-\tilde{A}_{w_j}^{-n}\hat{\bm{x}}^{(j)}_n\big)=A_{\tau_1}^{-1}\big(\tilde{A}_{w_1}^{-\ell_{n+1}}\bm{x}^{(1)}_{\ell_{n+1}}-\tilde{A}_{w_1}^{-\ell_n}\bm{x}^{(1)}_{\ell_n}\big).$$
It is easy to check that $\{\hat{\bm{x}}_n^{(j)}\}_{n\geq 0}\in l^2(\mathbb{C}^k,\tilde{A}_{w_j};r_{w_j}^{\sigma/2}\mu_{w_j}^{(\sigma-1)/2})$ and $\{\bm{x}_n^{(j)}-\hat{\bm{x}}_n^{(j)}\}_{n\geq 0}\in l^2(\mathbb{C}^k;r_{w_j}^{\sigma/2}\mu_{w_j}^{(\sigma-1)/2})$ using Lemma \ref{lemma35} if $\big(\{\bm{x}_n^{(1)}\}_{n\geq0},\{\bm{x}_n^{(2)}\}_{n\geq0},\cdots,\{\bm{x}_n^{(m)}\}_{n\geq0}\big)\in X_{p,\sigma}$ and $\sigma\neq 2l+\frac{d_S}{2},\forall l\geq 0$. The rest of the construction is easy and left to the reader.\vspace{0.1cm}

By applying the isomorphisms in Claim 1 and Claim 2, we finally are able to give a neat description of pretangents for $H^{2l}(K), {0\leq l\leq k}$. Theorem \ref{thm54} is  true since we can still show that the space of pretangents is stable under complex interpolation. {To show Theorem \ref{thm42}, a similar argument as Theorem \ref{thm410} is enough, noticing that we did not use (\textbf{A1}) in the proof of Lemma \ref{lemma49}.} The definition of $H_{00}^\sigma(K)$ remains the same even if (\textbf{A1}) is not satisfied, so Theorem \ref{thm56} remains the same. We may take $\tilde{H}_0^\sigma(K)$ as the right side of (\ref{eqnthm410}), then we can see that $H_{00}^\sigma(K)\subset H_0^\sigma(K)\subset \tilde{H}_0^\sigma(K)$, and $[\tilde{H}^0_{\sigma_1}(K),\tilde{H}^0_{\sigma_2}(K)]_\theta=H_{00}^\sigma(K)$ with $\sigma=(1-\theta)\sigma_1+\theta\sigma_2$. Theorem \ref{thm57} then follows as well {using Theorem \ref{thm56}}. Lastly, since Theorem \ref{thm62} is a consequence of the above theorems, it remains valid. 

\bibliographystyle{amsplain}

%\clearpage
%\newpage
\end{document}